\documentclass{amsart}

\usepackage[]{graphicx}
\usepackage[]{color}

\newcommand{\cl}{C \kern -0.1em \ell}

\newcommand{\BC}{{\rm \kern.24em \vrule width.05em height1.4ex depth-.05ex
            \kern-.26em C}}
\newcommand{\BZ}{{\rm \kern.24em \vrule width.05em height1.4ex depth-.05ex
            \kern-.26em Z}}
\newcommand{\BR}{{\rm\hskip 0.1pt I\hskip -2.15pt R}}
\newcommand{\BN}{{\rm\hskip 0.1pt I\hskip -2.15pt N}}

\newcommand{\vv}{{\bf v}}
\newcommand{\e}{{\bf e}}

\newtheorem{theorem}{Theorem}[section]
\newtheorem{example}{Example}[section]
\newtheorem{remark}{Remark}[section]
\newtheorem{lemma}{Lemma}[section]
\newtheorem{proposition}{Proposition}[section]
\newtheorem{corollary}{Corollary}[section]
\newtheorem{definition}{Definition}[section]

\begin{document}


\title[(Discrete) Almansi Type Decompositions]{(Discrete) Almansi Type Decompositions: \\ An umbral calculus framework based on $\mathfrak{osp}(1|2)$ symmetries}

\author{N.~Faustino}
\address{Centre for Mathematics,
University of Coimbra,
Largo D. Dinis,
Apartado 3008,
P-3001 - 454 Coimbra
, Portugal
}
\email{nelson@mat.uc.pt}
\thanks{N. Faustino was supported by {\it Funda\c c\~ao para a Ci\^encia e a Tecnologia} under the fellowship
SFRH/BPD/63521/2009 and the project PTDC/MAT/114394/2009.}

\author{G. Ren}
\address{University of Science and Technology of China,
 Department of Mathematics, Hefei, Anhui 230026, P. R. China}
\email{rengb@ustc.edu.cn}
\thanks{G. Ren was supported by the \textit{Wu Wen-Tsun Key Laboratory of Mathematics}, USTC, Chinese Academy of Sciences and by
the NNSF of China (No. 11071230).}

\date{\today}

\begin{abstract}
We introduce the umbral calculus formalism for hypercomplex variables starting from the fact that the algebra of multivariate polynomials $\BR[\underline{x}]$ shall be described in terms of the generators of the Weyl-Heisenberg algebra.
The extension of $\BR[\underline{x}]$ to the algebra of Clifford-valued polynomials $\mathcal{P}$ gives rise to an algebra of Clifford-valued operators whose canonical generators are isomorphic to the orthosymplectic Lie algebra $\mathfrak{osp}(1|2)$.

This extension provides an effective framework in continuity and discreteness that allow us to establish an alternative formulation of Almansi decomposition in Clifford analysis (cf. \cite{Ryan90,MR02,MAGU}) that corresponds to a meaningful generalization of Fischer decomposition for the subspaces $\ker (D')^k$.

We will discuss afterwards how the symmetries of $\mathfrak{sl}_2(\BR)$ (even part of $\mathfrak{osp}(1|2)$) are ubiquitous on the recent approach of \textsc{Render} (cf. \cite{Render08}), showing that they can be interpreted in terms of the method of separation of variables for the Hamiltonian operator in quantum mechanics.
\end{abstract}

\subjclass[2010]{30G35;35C10;39A12;81S05}

\keywords{Almansi theorem,~Clifford analysis,~hypercomplex variables,~orthosymplectic Lie algebras,~umbral calculus}


\maketitle

\section{Introduction}

\subsection{The Scope of Problems}
In the last two decades considerable attention has been given to the study of polynomial sequences for hypercomplex variables in different contexts.
For example, in the approach proposed by \textsc{Faustino~\& K\"ahler} (cf. \cite{FK07}), rising and lowering factorials yield e.g. the classical Bernoulli and Euler polynomials (cf. \cite{Tempesta08}) are the discrete analogue of homogeneous polynomials that appear in Fischer's decomposition involving difference Dirac operators. The hypercomplex generalization of this polynomials was studied recently in \cite{MT08} by \textsc{Malonek \& Tomaz} in connection with Pascal matrices.

Roughly speaking, the construction of hypercomplex Bernoulli polynomials shall be obtained via Appell sets \cite{CM08,MF08}. This was investigated at an early stage by \textsc{Abul-ez \& Constales} in \cite{AC90} in terms of basic sets of hypercomplex polynomials. In this case, the resulting hypercomplex polynomials shall be interpreted as the Cauchy-Kovaleskaya extension of the rising factorials considered in \cite{FK07}. For a fully explanation of Cauchy-Kovaleskaya extension we refer to \cite{DSS} (Subsection II.5); a meaningful characterization of Cauchy-Kovaleskaya's extension in interplay with Segal-Bargmann spaces can be found in \cite{CnopsKisil99} (Subsection 2.2).

Along with the construction of basic polynomial sequences in hypercomplex variables, other directions have been followed to construct Appell sets, namely by taking the Fueter-Sce extension of complex monomials $z^k$ (cf. \cite{Norman09}), using a Fourier series expansion of square integrable monogenic functions on the unit ball (cf. \cite{BG10}) or alternatively using a Gelfand-Tsetlin basis approach (cf. \cite{BGLS10}) that essentially is a combination between Fischer decomposition and Cauchy-Kovalevskaya extension.

 With respect to the discrete setting, it was further developed in the Ph.D thesis of \textsc{Faustino} (cf. \cite{faustino09}) that families of discrete polynomials can be constructed as a
blending between continuous and discrete Clifford analysis giving an affirmative answer to the paper of \textsc{Malonek \& Falc\~ao} (cf. \cite{MF08}).

According to this proposal, discrete Clifford operators underlying the orthogonal group $O(n)$ were introduced by means of representations of the Lie superalgebra $\mathfrak{osp}(1|2)$. Moreover, the refinement of discrete harmonic analysis follows from the representation of the Lie algebra $\mathfrak{sl}_2(\BR)$ as the even part of $\mathfrak{osp}(1|2)$ while the blending between {\it continuum} and discrete Clifford analysis was obtained via a Sheffer map (cf.~\cite{RR,roman84}) that essentially maps the homogeneous polynomials onto basic polynomial sequences of binomial type.

This approach combines the radial algebra based approach proposed by \textsc{Sommen} in \cite{Sommen97} with the umbral calculus approach postponed in \cite{BLR98} by \textsc{Di Bucchianico,Loeb \& Rota}. The main novelty of this new approach rests mostly from the fact that continuous and discrete Clifford analysis are described as realizations of the well-known Wigner quantum systems (cf. \cite{Wigner50}) on which the Sheffer map shall be interpreted as a gauge transformation that keeps invariant the symmetries of both systems (see e.g. \cite{faustino10} for a sketch of this approach).

There were still to consider alternative constructions of discrete Clifford analysis using different type of symmetries. One of them proposed in the preprint of \textsc{Faustino \& K\"ahler} (cf.~\cite{FK08}), Clifford analysis on symmetric lattices corresponds to a mimetic description of Hermitian Clifford analysis on which the unitary group $U(n)$ appears as the natural candidate to the induced representations for the algebra of Clifford-valued operators (cf. \cite{BSES10}). The major obstacle arising this construction follows from the fact that the multiplication operators $x_jT_{h}^{-j}$ and $x_jT_{h}^{+j}$ do not commute and hence there is no chance to get a radial algebra structure (cf. \cite{Sommen97}). For a complete survey besides this drawback we refer to \cite{faustino09} (Section 3).

 Quite recently, in the recent approach of \textsc{De Ridder, De Schepper, K\"ahler \& Sommen} (cf.~\cite{RSKS10}) the Weyl-Heisenberg symmetries encoded in the forward/backward finite difference operators $\partial_h^{\pm j}$ and multiplication operators $x_jT_{h}^{\mp j}$ were replaced by `skew'-Weyl symmetries with the purpose to get, in analogy with the Hermitian setting, linear independence between the vector multiplication operators $X^+=\sum_{j=1}^n\e_j^+ X_j^+$ and $X^-=\sum_{j=1}^n\e_j^- X_j^-$. As a result, the authors have shown that the Euler polynomials are the resulting discrete polynomials that yield a Fischer decomposition for $D_h^{+}+D_h^{-}=\sum^n_{j=1} \e_j^+ \partial_h^{+j}+\e_j^- \partial_h^{-j}$. Besides this approach there is an open question regarding the group of induced representations of such system.

Let us turn now our attention for the Almansi decomposition state of art in Clifford and harmonic analysis. The theorem formulated below:
\medskip

\textbf{Almansi's Theorem} (cf.~\cite{ALM,ACL}) {\em {If $f$ is
polyharmonic of degree $k$ in a starlike domain with center $0$,
then there exist uniquely defined harmonic functions on $\Omega$ $f_0, \cdots,f_{k-1}$ such that
$$f(x)=f_0(\underline{x})+|x|^2 f_2(\underline{x})+\cdots+ |x|^{2(k-1)}f_{k-1}(\underline{x}).$$}}
corresponds to the Almansi decomposition for polyharmonic functions.

One can find important applications and generalizations of this
result for several complex variables in the monograph of
\textsc{Aronszajn, Creese \& Lipkin,} \cite{ACL}, e.g. concerning
functions holomorphic in the neighborhood of the origin in $\BC^n$.

In the harmonic analysis setting, the importance of this result was recently explored by \textsc{Render} in \cite{Render08}, showing that for functions belonging to the real Bargmann space, there is an intriguing connection between the existence of a Fischer inner pair (cf. \cite{Render08}) with the problem of uniqueness for polyharmonic functions posed by \textsc{Hayman} in \cite{Hay94} (cf. \cite{Render08},~Section 9) and as a characterization for the entire solutions of Dirichlet's problem (cf. \cite{Render08},~Section 10).

  In the Clifford analysis setting, the Almansi theorem shall be understood as a meaningful generalization of Fischer decomposition for hypercomplex variables without requiring \textit{a-priori} a Fischer inner product (cf. \cite{DSS}, pp. 204-207).
   This result plays a central role in the study of polymonogenic functions likewise in the study of polyharmonic functions as refinements of polymonogenic functions. This was consider in the begining of $90's$ by \textsc{Ryan} \cite{Ryan90} to study the invariance of iterated Dirac operators in relation to M\"obius transformations on manifolds. On the last decade \textsc{Malonek \& Ren} established a general framework which describe the decomposition of iterated kernels for different function classes \cite{MR02,MAGU}. Besides the approach of \textsc{Cohen, Colonna, Gowrisankaran~\& Singman} (cf.~\cite{CCGS}) regarding polyharmonic functions on trees and the approaches on Fischer decomposition for difference Dirac operators proposed by \textsc{Faustino \& K\"ahler} \cite{FK07} and \textsc{De Ridder, De Schepper, K\"ahler \& Sommen} \cite{RSKS10}, up to now there is no established framework on Almansi type decompositions as a general method for obtaining special representations for discrete hypercomplex functions.

\subsection{Motivation of this approach}
The umbral calculus formalism proposed by \textsc{Roman \& Rota} (cf. \cite{roman84,RR}) have received in the last fifteen years the attention of mathematicians and physicists. Besides the papers of \textsc{Di Bucchianico \& Loeb} (cf.~\cite{BL96}) and \textsc{Di Bucchianico,~Loeb, \& Rota} (cf.~\cite{BLR98}) devoted to classical aspects of umbral calculus, further applications were developed after the papers of \textsc{Smirnov \& Turbiner} (cf. \cite{SmirnovTurbiner95}) and \textsc{Dimakis, M\"uller-Hoissen \& Striker} in the mid of the $90's$ (cf. \cite{DHS96}) with special emphasis to the systematic discretization of Hamiltonian operators preserving Weyl-Heisenberg symmetries (cf. \cite{ML01,LTW04}), to the construction of Appell sets (cf. \cite{Tempesta08}) and complete orthogonal systems of polynomials (cf. \cite{DLW08}) based on the theory of Sheffer sets likewise to the solution of the Boson-Normal ordering problem in quantum mechanics based on the interplay between combinatorial identities using binomial sums with the construction of coherent states (cf. \cite{BlasiakDattoliHorzelaPenson06}).

When we take the tensor product between the algebra of multivariate polynomials $\BR[\underline{x}]$ with the Clifford algebra of signature $(0,n)$ in $\BR^n$, the resulting algebra of Clifford-valued polynomials is described in terms of Lie symmetries underlying the Lie algebra $\mathfrak{sl}_2(\BR)$ and the Lie superalgebra $\mathfrak{osp}(1|2)$  (see \cite{DBieSommen07,faustino10,CFK10} and the references given there) while the Fischer decomposition of the algebra of homogeneous Clifford-valued polynomials in terms of spherical harmonics and spherical monogenics follows from the Howe dual pair technique (see \cite{BGLS10} and references given there) applied to $\mathfrak{sl}_2(\BR)\times O(n)$ and $\mathfrak{osp}(1|2)\times O(n)$, respectively. The details of such technique can also be found in the papers of \textsc{Howe} (cf.~\cite{Howe89}) and \textsc{Cheng \& Zhang} (cf.\cite{ChengZhang04}).

Although in the last years the Lie (super)algebra framework was successfully applied in Clifford analysis, these kinds of algebras are also ubiquitous e.g. in the old works of \textsc{Wigner} (cf. \cite{Wigner50}) and \textsc{Turbiner} (cf. \cite{Turbiner88}):

In \cite{Wigner50} it was shown that the so-called Wigner quantum systems that describes a motion of a particle on the ambient space $\BR^n$ may be characterized in terms of symmetries of $\mathfrak{osp}(1|2n)$; in \cite{Turbiner88} the eigenfunctions for the Hamiltonian operators were computed explicitly by taking into account the $\mathfrak{sl}_2(\BR)$ symmetries of such system while the eigenvalues were described as an infinite number of (unplaited) sheets lying on a Riemann surface. Quite recently, in \cite{Zhang08} \textsc{Zhang} also apply this framework to the study of quantum analogues for the Kepler problem in the superspace setting.

\subsection{Organization of the paper}
In this paper we will derive an umbral counterpart for the well known Almansi type decomposition for hypercomplex variables by employing combinatorial and algebraic techniques regarding umbral calculus (cf.~\cite{DHS96,BLR98}), radial algebras (cf. \cite{Sommen97}) and the Howe dual pair technique confining nonharmonic analysis and quantum physics (cf. \cite{Howe85,HT92}).

We will start to introduce the umbral calculus framework in the algebra of Clifford-valued polynomials $\mathcal{P}:=\BR[\underline{x}]\otimes \cl_{0,n}$ as well as the symmetries preserved under the action of the Sheffer map, showing that there is a mimetic transcription of classical Clifford analysis to discrete setting that generalizes complex analysis to higher dimensions (cf. \cite{DSS,GM}).
Roughly speaking, in umbral calculus the algebra of polynomials $\BR[\underline{x}]$ can be recognized as being isomorphic
to the algebra generated by position and momentum operators $x_j'$ and $O_{x_j}$, respectively, satisfying the Weyl-Heisenberg relations
\begin{eqnarray}
\label{weylh}[O_{x_j},O_{x_k}]=0=[x_j',x_k'], & [O_{x_j},x_k']=\delta_{jk}{\bf id}.
\end{eqnarray}
Here and elsewhere $[{\bf a},{\bf b}]:={\bf a}{\bf b}-{\bf b}{\bf a}$ denotes the commuting bracket between ${\bf a}$ and ${\bf b}$.

Moreover, if we take
$\e_1,\ldots, \e_n$ as the Clifford algebra generators satisfying the anti-commuting
relations
$\{\e_j,\e_k\}:=\e_j\e_k+\e_k\e_j=-2\delta_{jk}$, umbral Clifford analysis (cf.~\cite{faustino10})
deals with the study of the algebra of
differential operators
\begin{eqnarray*}
\mbox{Alg}\left\{
x_j',O_{x_j},\e_j~:~j=1,\ldots,n\right\},
\end{eqnarray*}
For a complete survey besides this approach we refer to \cite{faustino09} (Section 3).

Moreover, taking into account the action of the orthosymplectic Lie algebra of type $\mathfrak{osp}(1|2)$ on the subspaces $(x')^s\ker D'$ we will derive some recursive relations (Lemma \ref{loweringLemma} and Proposition \ref{mappingProperty}) and inversion formulae (Lemmata \ref{Is'Es'} and \ref{le:51}) which allow us to decompose the subspace $\ker (D')^k$ (the so-called umbral polymonogenic functions of degree $k$) as a direct sum of subspaces of the type $P_s(x')\ker D'$, for $s=0,1,\ldots,k-1$, where $P_s(x')$ stands a polynomial type operator of degree $s$ satisfying the mapping property $P_s(x'): \ker (D')^s \rightarrow \ker D'$.

This in turn gives an alternative interpretation for the results obtained by \textsc{Ryan} (cf.~\cite{Ryan90}), \textsc{Malonek \& Ren} (\cite{MR02,MAGU}) and \textsc{Faustino \& K\"ahler} \cite{FK07} in terms of the symmetries of $\mathfrak{osp}(1|2)$.

Finally, in Subsection \ref{HarmonicOscillator} we will give an interpretation for the recent approach of \textsc{Render} (cf. \cite{Render08}) showing that for a special choice of the potential operator $V_\hbar(x')$ the Almansi decomposition encoded in the quantized Fischer inner pair $(\left(2V_\hbar(x')\right)^k,(\Delta')^k)$ is nothing else than a $\mathfrak{sl}_2(\BR)$ based diagonalization of the Hamiltonian $\mathcal{H}'=-\frac{1}{2}\Delta'+V_\hbar(x')$.

\section{Umbral Clifford Analysis}\label{umbralCliffordAnalysis}

\subsection{Umbral calculus revisited}\label{umbralC}
In this section we will review some basic notions regarding umbral calculus. The proof of further results that we will omit can be found e.g.~in \cite{RR,roman84,BL96} or alternatively in \cite{faustino09} (see Chapter 1).

In the following, we will set $\BR[\underline{x}]$ as the ring of polynomials over $\underline{x}=(x_1,x_2,\ldots,x_n)\in \BR^n$,
by $\alpha=(\alpha_1,\alpha_2,\ldots \alpha_n)$ the multi-index over $\BN^n_0$ and by $\underline{x}^{\alpha}=x_1^{\alpha_1}x_2^{\alpha_2}\ldots
x_n^{\alpha_n}$ the multivariate monomial(s) over
$\underline{x}$.
The partial derivative with respect to $x_j$ will be denoted by
$\partial_{x_j}:=\frac{\partial}{\partial{x_j}}$ while the $n-$tuple $\partial_{\underline{x}}:=(\partial_{x_1},\partial_{x_2},\ldots,\partial_{x_n})$ corresponds to the so-called gradient operator.

Here and elsewhere, we will also consider the following
notations:
\begin{eqnarray*}
\begin{array}{ccccc}
 \partial^{\alpha}_{\underline{x}}:=\partial_{x_1}^{\alpha_1}\partial_{x_2}^{\alpha_2}\ldots\partial_{x_n}^{\alpha_n}, & \alpha!=\alpha_1! \alpha_2 !\ldots \alpha_n! , &
 \left(\begin{array}{cc} \beta \\ \alpha \end{array}\right)=\frac{\beta!}{\alpha ! (\beta-\alpha)!}, &
 |\alpha|=\sum_{j=1}^n\alpha_j.
 \end{array}
\end{eqnarray*}

In addition we will denote by $\mbox{End}(\BR[\underline{x}])$ the algebra of linear operators acting on $\BR[\underline{x}]$.

By means of the differentiation formulae $\partial^{\alpha}_{\underline{x}}\underline{x}^{\beta}=0$ for $|\alpha|>|\beta|$ and
 $\partial^{\alpha}_{\underline{x}}\underline{x}^{\beta}=\frac{\beta!}{(\beta-\alpha)!}~\underline{x}^{\beta-\alpha}$ for $|\alpha|\leq |\beta|$, it turns out the representation of the binomial formula in terms of $\partial_{\underline{x}}$:
\begin{eqnarray}
\label{binomialFormula}(\underline{x}+\underline{y})^{\beta}=\sum_{|\alpha|=0}^{|\beta|}\left(\begin{array}{cc} \beta \\ \alpha \end{array}\right)\underline{x}^{\alpha}\underline{y}^{\beta-\alpha}
=\sum_{|\alpha|=0}^{\infty}\frac{[\partial_{\underline{x}}^\alpha\underline{x}^\beta]_{\underline{x}
    =\underline{y}}}{\alpha !}\underline{x}^{\alpha}
\end{eqnarray}

Linearity arguments shows that the extension of the above formula to $\BR[\underline{x}]$ is given by $f(\underline{x}+\underline{y})=\exp(\underline{y} \cdot
\partial_{\underline{x}})f(\underline{x})$, where $\exp(\underline{y} \cdot
\partial_{\underline{x}})=\sum_{|\alpha|=0}^{\infty} \frac{\underline{y}^\alpha}{\alpha!}\partial_{\underline{x}}^\alpha$ denotes the formal power series representation for the shift operator $T_{\underline{y}}f(\underline{x})=f(\underline{x}+\underline{y})$.

An operator $Q \in \mbox{End}(\BR[\underline{x}])$ is shift-invariant if and only if it commutes with  $T_{\underline{y}}=\exp(\underline{y} \cdot
\partial_{\underline{x}})$ for all $P \in \BR[\underline{x}]$ and $\underline{y}\in
 \BR^n$: $$[Q,T_{\underline{y}}]P(\underline{x}):=Q(T_{\underline{y}}P(\underline{x}))-T_{\underline{y}}(Q~P(\underline{x}))=0.$$

Under the shift-invariance condition for $Q$, the first expansion theorem (cf. \cite{BL96}) states that
any linear operator $Q: \BR[\underline{x}]\rightarrow
\BR[\underline{x}]$ is shift-invariant if and only if $Q$ is given in terms of the following formal power series expansion:
\begin{eqnarray*}
Q=\sum_{|\alpha|=0}^{\infty}\frac{a_\alpha}{\alpha!} \partial_{\underline{x}}^{\alpha},& \mbox{with}~~a_\alpha=\left[{Q \underline{x}^\alpha}\right]_{\underline{x}=\underline{0}}.
\end{eqnarray*}

Set $O_{\underline{x}}=(O_{x_1},O_{x_2},\ldots,O_{x_n})$ as a multivariate operator. We say that $O_{\underline{x}}$ is shift-invariant if and only if $O_{x_1},O_{x_2},\ldots,O_{x_n}$ are shift-invariant too.
Moreover, $O_{\underline{x}}$ is a multivariate delta operator if and only if there is a non-vanishing constant $c$ such that $O_{x_j}(x_k)=c \delta_{jk}$, holds for all $j,k=1,2,\ldots,n$.
It can be shown that if $O_{\underline{x}}$ is a multivariate delta operator, each $O_{x_j}$ lowers the degree of $P(\underline{x}) \in \mathbb{R}[\underline{x}]$. In particular $O_{x_j}(c)=0$ for each non-vanishing constant $c$ (cf. \cite{faustino09}, Lemmata 1.1.8 and 1.1.9) and hence, any multivariate delta operator $O_{\underline{x}}$
uniquely determines a polynomial sequence of binomial type $\{V_\alpha(\underline{x})~:~\alpha \in \mathbb{N}_0^n\}$,
 (cf.~\cite{faustino09}~Theorems 1.1.12 and 1.1.13):
 \begin{eqnarray}
\label{binomialType}V_\beta(\underline{x}+\underline{y})=\sum_{|\beta|=0}^{|\alpha|}\left(\begin{array}{cc} \beta \\ \alpha \end{array}\right)V_\alpha(\underline{x})V_{\beta-\alpha}(\underline{y})
\end{eqnarray}
such that $V_{\underline{0}}(\underline{x})=1$, $V_{\alpha}(\underline{0})=\delta_{\alpha,0}$ and $O_{x_j}V_{\alpha}(\underline{x})=\alpha_jV_{\alpha-\vv_j}(\underline{x}),$
where $\vv_j$ stands the $j-$element of the canonical basis of $\mathbb{R}^n$.

The \index{Pincherle derivative}Pincherle derivative of $O_{x_j}$ with respect to $x_j$
is defined formally as the commutator between $O_{x_j}$ and $x_j$:
\begin{eqnarray*}
\label{Pincherle} O'_{x_j}f(\underline{x}):=[O_{x_j},x_j]f(\underline{x})=O_{x_j}(x_jf(\underline{x}))-x_j(O_{x_j}f(\underline{x})).
\end{eqnarray*}
This canonical operator plays an important role in the construction of basic polynomial sequences of binomial type (cf. \cite{DHS96,DLW08,BlasiakDattoliHorzelaPenson06}).
The subsequent results allows us to determine in which conditions $(O_{x_j}')^{-1}$ exists.

We will start with the following lemma:

\begin{lemma}\label{shiftinvPincherle}
The \index{Pincherle derivative}Pincherle derivative of a shift-invariant operator\index{shift-invariant operator} $Q$ is shift-invariant.
\end{lemma}

Regardless the last lemma one looks to shift-invariant operators $Q$
as formal power series $Q(\underline{x})=\sum_{\alpha}\frac{a_\alpha}{\alpha!}{\underline{x}}^\alpha$ obtained {\it viz} the replacement of $\underline{x}$ by $\partial_{\underline{x}}$, i.e. $\iota[Q(\underline{x})]=Q(\partial_{\underline{x}})$ where $\iota:\widehat{\BR[\underline{x}]} \rightarrow \mbox{End}(\BR[\underline{x}])$ is defined as a mapping between the algebra of formal power series $\widehat{\BR[\underline{x}]}$ and the algebra of linear operators acting on $\BR[\underline{x}]$.
According to the isomorphism theorem (see \cite{roman84},~Theorem 2.1.1.), $\iota$ is one-to-one and onto. This in turn leads to the following proposition:

\begin{proposition}\label{inv}
A shift-invariant operator\index{shift-invariant operator} $Q$ has its inverse if and only if $Q 1 \neq 0$.
\end{proposition}

From Lemma \ref{shiftinvPincherle} notice that $O_{x_j}'$ is shift-invariant whenever $O_{x_j}$ is shift-invariant. Since from the definition $O_{x_j}'({\bf 1})=O_{x_j}(x_j)$ and $O_{x_j}(x_j)$ is a non-vanishing constant, Proposition \ref{inv} asserts that $(O_{x_j}')^{-1}$ exists locally as a formal series expansion involving multi-index derivatives $\partial_{\underline{x}}^\alpha$.

The former description in terms of Pincherle derivatives allows us to determine $V_\alpha(\underline{x})$ as a polynomial sequence obtained from the action of $(\underline{x}')^{\alpha}:=\prod_{k=1}^n(x_k')^{\alpha_k}$, with $x_k':=x_k (O_{x_k}')^{-1}$, on the constant polynomial $\Phi={\bf 1}$ (see also \cite{roman84}, page 51, Corollary 3.8.2):
\begin{eqnarray}
\label{rodrigues}V_\alpha(\underline{x})=(\underline{x}')^{\alpha}{\bf
1}.
\end{eqnarray}

The properties of basic polynomial sequences are naturally characterized
within the extension of the mapping property $\Psi_{\underline{x}}: \underline{x}^\alpha \mapsto
V_\alpha(\underline{x})$ to $\mathbb{R}[\underline{x}]$. According to \cite{roman84}, this mapping is the well-known Sheffer map that link two basic polynomial sequences of binomial type. It is clear from
the construction that
$\Psi_{\underline{x}}^{-1}$ exists and it is given by the linear extension of
$\Psi_{\underline{x}}^{-1}: V_\alpha(\underline{x}) \mapsto
\underline{x}^\alpha$ to $\mathbb{R}[\underline{x}]$. In addition we get the following properties on $\mathbb{R}[\underline{x}]$:
\begin{eqnarray*}
\label{intertwining}O_{x_j}
=\Psi_{\underline{x}} \partial_{x_j}\Psi_{\underline{x}}^{-1}, &~~ \mbox{and}~~ & x_j'
=\Psi_{\underline{x}} x_j \Psi_{\underline{x}}^{-1}.
\end{eqnarray*}

From the border view of quantum mechanics, the $2n+1$ operators $x_1',\ldots,x_n'$,
$O_{x_1},\ldots,O_{x_n}$ and ${\bf id}$ generate the \index{Bose algebra}Bose algebra isomorphic to $\BR[\underline{x}]$ (cf. \cite{BLR98}). Indeed
for $\Phi={\bf 1}$ (the so-called vacuum vector)
$O_{x_j}(\Phi)=0$ holds for each $j=1,\ldots,n$ while the raising and lowering operators, $x_j':V_{\alpha}(\underline{x}) \mapsto V_{\alpha+\vv_j}(\underline{x})$ and $O_{x_j}:V_{\alpha}(\underline{x}) \mapsto \alpha_j~V_{\alpha-\vv_j}(\underline{x})$ respectively, satisfying the Weyl-Heisenberg relations given by (\ref{weylh}).

Due to this correspondence, we would like to stress that the quantum mechanical description of umbral calculus give us many degrees of freedom to construct the raising operators $x_j':V_{\alpha}(\underline{x}) \mapsto
V_{\alpha+\vv_j}(\underline{x})$ in such way that the commuting relations (\ref{weylh}) fulfil.
 In particular, in \cite{DHS96,DLW08} it was pointed out the importance to
consider the following symmetrized versions of $x_j (O_{x_j}')^{-1}$
\begin{eqnarray}
\label{xj'2} x_j'=\frac{1}{2}(x_j (O_{x_j}')^{-1}+(O_{x_j}')^{-1}x_j)
\end{eqnarray}
as a special type of canonical discretization.

\subsection{Basic operators}\label{umbralCliffordOperators}

In what follows we will use the notation introduced in Section
\ref{umbralC}. In addition 
we introduce $\cl_{0,n}$ as the
algebra determined by the set of vectors $\e_1,\e_2,\ldots,\e_n$
satisfying the graded relations with respect to the anti-commuting bracket $\{{\bf a},{\bf b}\}={\bf a}{\bf b}+{\bf b}{\bf a}$:
\begin{eqnarray}
\label{CliffPinch}\{
\e_j,\e_k\}=-2\delta_{jk},
\end{eqnarray}

The above algebra is commonly known in literature as the Clifford algebra of
signature $(0,n)$ (cf. \cite{DSS},~Chapters 0 \& I; \cite{GM}, Chapter 1) which corresponds to a particular example of an algebra of radial type. Indeed the anti-commutator $\{
\e_j,\e_k\}$ (scalar-valued quantity) commutes with all the basic vectors $\e_j$:
\begin{eqnarray}
\label{radial}\left\{~\left[\e_j,\e_k \right],\e_l~\right\}=0, &
\mbox{for all}~j,k,l=1,\ldots,n.
\end{eqnarray}
For further details concerning the construction of $\cl_{0,n}$ as an algebra of radial type we refer to \cite{Sommen97}.

 Additionally, we will denote by $\mathcal{P}=\BR[\underline{x}]\otimes \cl_{0,n}$ the algebra of Clifford-valued polynomials and by $\mbox{End}(\mathcal{P})$ the algebra of linear operators acting on $\mathcal{P}$.
The Weyl-Heisenberg character of the operators $x_j'$ and $O_{x_j}$ combined with the radial character of the generators underlying $\cl_{0,n}$ allows us to define
umbral Clifford analysis
as the study of the following algebra of
differential operators:
\begin{eqnarray}\label{eq:1.43}
\mbox{Alg}\left\{
x_j',O_{x_j},\e_j~:~j=1,\ldots,n\right\}
\end{eqnarray}
Furthermore, the umbral counterparts for the Dirac operator, vector variable and Euler operator, $D'$, $x'$ and $E'$ respectively, defined as follows
\begin{eqnarray}
\label{dirac}D'=\sum_{j=1}^n \e_j O_{x_j}, \\
\label{x'}x'=\sum_{j=1}^n \e_j x_j', \\
\label{E'} E'=\sum_{j=1}^n x_j' O_{x_j},
\end{eqnarray}
correspond to linear combinations of the elements of the algebra defined in (\ref{eq:1.43}).

In this context, the operators (\ref{dirac})-(\ref{E'}) shall be understood as basic left endomorphisms acting on the algebra $\mbox{End}(\mathcal{P})$.
Along this paper we will use several times the notation $\Delta':=\sum_{j=1}^n  O_{x_j}^2$ to refer the umbral counterpart of the Laplace operator $\Delta=\sum_{j=1}^n \partial_{x_j}^2$.

The next lemma naturally follows from straightforward computations obtained by direct combination of relations (\ref{weylh}) and (\ref{CliffPinch}):
\begin{lemma}[cf. \cite{faustino09}, Lemma 3.4.3, pp. 68]\label{x&D}
The operators $x',D',E' \in \mbox{End}(\mathcal{P})$ satisfy the following anti-commutation relations
\begin{eqnarray*}
\{ x',x'\}=-2\sum_{j=1}^n (x_j')^2, &\{D',D'\}=-2\Delta', & \{x',D'\}=-2E'-n{\bf id}.
\end{eqnarray*}
\end{lemma}

From the first relation of Lemma \ref{x&D}, $-(x')^2({\bf 1})$ is a scalar-valued quantity while from the second relation $\Delta':=-(D')^2$ is a second order operator satisfying the vanishing condition $\Delta'({\bf 1})=0$.
On the other hand, in the third relation the action $E'$ on $\mbox{End}(\mathcal{P})$ can be rewritten as a the following identity involving the anti-commuting relation between $x'$ and $D'$:
\begin{eqnarray}
\label{E'&x'D'}E'=\sum_{j=1}^n x_j' O_{x_j}=-\frac{1}{2}\left(\{x',D'\}+n{\bf id}\right).
\end{eqnarray}

The action on both sides of the above identity on $\Phi={\bf 1}$ allows us to recast the dimension of the ambient space $\BR^n$ as $n=-D'(x'{\bf 1})$. So, the polynomial $\Phi={\bf 1}$ shall be interpreted as the corresponding ground level eigenstate while the dimension of the ambient space $\BR^n$ appears as twice of the
ground level energy associated to the harmonic oscillator containing
$n$ degrees of freedom.

It is also clear from Lemma \ref{x&D} that the anti-commutators $\{ x',x'\}$, $\{ D',D'\}$ and $\{ x',D'\}$ are scalar-valued.
Thus, from (\ref{dirac}),(\ref{x'}) and (\ref{E'&x'D'}) the operators $(x')^2$,~$-(D')^2$ and $E'$ shall be view as generalizations for the
norm squared of a vector variable in the Euclidean space, the Laplacian operator and Euler operator, respectively.
From the border view of quantum mechanics, the
operators $-\frac{1}{2}(x')^2$ and $-\frac{1}{2}\Delta'$ describe a
spherical potential and the kinetic energy, respectively.

Here we would like also to stress that the operator $E'$ (see identity (\ref{E'})) comprises at the same time the concept of directional derivative introduced by \textsc{Howe} (cf.~\cite{Howe85}) for quantum groups with the concept of non-shift-invariant mixed/number operator given by \textsc{Di Bucchianico, Loeb \& Rota} (cf. \cite{BLR98}).

We will end this subsection by exploring and discussing some examples regarding the construction of the operators (\ref{dirac})-(\ref{E'}).

\begin{example}\label{DiracOperators}
If we take $O_{x_j}=\partial_{x_j}$, $D'$ and $x'$
 coincide with the standard Dirac and coordinate variable operators, respectively:
 \begin{eqnarray*}
 D=\sum_{j=1}^n \e_j \partial_{x_j}, & x=\sum_{j=1}^n \e_j x_j.
 \end{eqnarray*}
while $E'$ corresponds to the classical Euler operator $E=\sum_{j=1}^nx_j \partial_{x_j}$.

    Furthermore, the {\it continuum} Hamiltonian $\frac{1}{2}\left(-\Delta+|\underline{x}|^2 \right)$
     can we rewritten as $$\frac{1}{2}\left(-\Delta+|\underline{x}|^2 \right)=\frac{1}{2}(D^2-x^2).$$
\end{example}

\begin{example}\label{differenceDiracOperators}
Next we will consider $D'$ as a difference Dirac operator given in terms of the forward differences $\partial_{h}^{+j}f(\underline{x})=\frac{f(\underline{x}+h{\bf v}_j)-f(\underline{x})}{h}$ supported on the grid $h\BZ$:
$$D'=\sum_{j=1}^n \e_j \partial_{h}^{+j}=\sum_{j=1}^n \e_j \frac{T_{h \vv_j}-{\bf id}}{h}.$$

The square of $D'$ corresponds to $(D')^{2}=-\frac{1}{h^2}\sum_{j=1}^n (T_{2h\vv_j}-2T_{h\vv_j}+{\bf id})$. On the other hand, the formal series expansion for $\partial_{h}^{+j}$ is given by $\partial_h^{+j}=\frac{1}{h}\left(\exp(h\partial_{x_j})-{\bf id}\right),$ and
$[\partial_h^{+j},x_j]=T_{h {\bf v}_j}=\exp(h\partial_{x_j})$ corresponds to the Pincherle derivative for $\partial_{h}^{+j}$.

Thus the operator $x'$ corresponds to
$$ x'=\sum_{j=1}^n \e_j~ x_j T_{-h {\bf v}_j}=\sum_{j=1}^n \e_j~ x_j \exp(-h\partial_{x_j}).$$

Alternatively, using relation (\ref{xj'2}), the operator $x'$ can also be taken as
$$ x'=\frac{1}{2}\sum_{j=1}^n \e_j ~\left(x_j T_{-h {\bf v}_j}+T_{-h {\bf v}_j}x_j \right)=\frac{1}{2}\sum_{j=1}^n \e_j ~\left(x_j \exp(-h\partial_{x_j})+\exp(-h\partial_{x_j})x_j \right).$$
Then we can consider two different constructions for $E'$ (see relations (\ref{E'}) and (\ref{E'&x'D'})):
\begin{itemize}
\item $E'=\sum_{j=1}^n x_j T_{-h \vv_j} \partial_h^{+j}=\sum_{j=1}^n x_j \partial_h^{-j}$;
\item $E'=\frac{1}{2}\sum_{j=1}^n \left(x_j T_{-h \vv_j}+T_{-h \vv_j}x_j\right)\partial_h^{+j}=\frac{1}{2}\sum_{j=1}^n x_j \partial_h^{-j}+\frac{1}{2}\sum_{j=1}^n T_{-h \vv_j }x_j \partial_h^{+j}$.
\end{itemize}
Hereby $\partial_h^{-j}=\frac{1}{h}({\bf id}-T_{-\vv_j})$ corresponds to the backward finite difference operator acting on the grid $h\BZ$.
\end{example}

\begin{example}\label{discreteClifford2}
Now we will replace the forward finite differences $\partial_h^{+j}$ used to define $D'$ in Example \ref{differenceDiracOperators} by a central difference operator acting on $h\BZ^n$:
$$O_{x_j}f(\underline{x})=\frac{f(\underline{x}+h\vv_j)-f(\underline{x}-h\vv_j)}{2h}.$$ The formal series expansion for these operators  correspond to $$O_{x_j}=\frac{1}{2h}\left(\exp\left(h\partial_{x_j}\right)-\exp\left(-h\partial_{x_j}\right)\right)=\frac{1}{2h}\sinh\left(h\partial_{x_j}\right)$$ and moreover the formal series expansion for $D'$ is given by
    \begin{eqnarray*}
    \label{DiscreteDirac}D'=\frac{1}{2h}\sum_{j=1}^n \e_j~ \sinh\left(h\partial_{x_j}\right).
    \end{eqnarray*}
    The square of $(D')^2$ splits the star laplacian on a equidistant grid with mesh-width $2h$:
    \begin{eqnarray*}
    \label{starLapl}-(D')^2=\sum_{j=1}^n
\frac{\exp(2h\partial_{x_j})-2{\bf id}+\exp(-2h\partial_{x_j})}{4h^2}=\sum_{j=1}^n
\frac{T_{2h \vv_j}-2{\bf id}+T_{-2h \vv_j}}{4h^2}.
 \end{eqnarray*}

Therefore, the construction of $x'$ and $E'$ shall be take into account the following formal series expansion for $O_{x_j}':$ $$O_{x_j}'f(\underline{x})=\frac{f\left(\underline{x}+h{\bf v}_j\right)+f\left(\underline{x}-h{\bf v}_j\right)}{2}=\cosh\left(h\partial_{x_j}\right)f(\underline{x}).$$

Using the relation $\cosh\left(h\partial_{x_j}\right)=\frac{1}{h}\exp\left(-h\partial_{x_j}\right)\left( {\bf id}-\exp\left(2h\partial_{x_j}\right) \right)$ combined with the standard Von Neumann series expansion of $\left( {\bf id}-\exp\left(2h\partial_{x_j}\right)  \right)^{-1}$, we get the following asymptotic expansion for $(O_{x_j}')^{-1}:$

$$ (O_{x_j}')^{-1}=-h\left( {\bf id}-\exp\left(2h\partial_{x_j}\right) \right)^{-1}\exp\left(h\partial_{x_j}\right)
=-h\sum_{k=0}^\infty \exp\left((2k+1)h\partial_{x_j}\right),
$$
or equivalently $(O_{x_j}')^{-1}=-h\sum_{k=0}^\infty T_{(2k+1)h \vv_j}$.
The above inverse only exists whenever $\left\|T_{2h \vv_j}\right\|=\left\|\exp\left(2h\partial_{x_j}\right)\right\|<1$.

Alternatively, we can express $(O_{x_j}')^{-1}$ using the following formal integral representation in terms of the Laplace transform $(\mathcal{L} f)(s)=\int_{0}^\infty e^{-st}f(t)~dt$ (cf. \cite{BLR98}):
$$ (O_{x_j}')^{-1}=-h\int_{0}^\infty e^{-st}\exp(h(2t+1)\partial_{x_j})~dt=-h\int_{0}^\infty e^{-st}T_{h(2t+1)\vv_j}~dt.$$
\end{example}

The umbral Dirac operator introduced in Example \ref{differenceDiracOperators} corresponds to the forward difference Dirac operator introduced by \textsc{Faustino \& K\"ahler} in \cite{FK07}.
Here we would like to notice that in Example \ref{differenceDiracOperators}, the square $(D')^2$ does not split the star Laplacian
\begin{eqnarray*}
\Delta_h f(\underline{x})= \sum_{j=1}^n \frac{f(\underline{x}+h{\bf v}_j)+f(\underline{x}-h{\bf v}_j)-2f(\underline{x})}{h^2},
\end{eqnarray*}
which means that discrete harmonic analysis can not be refined in terms of discrete Dirac operators involving only forward differences (cf. \cite{FK07}).

As we see in Example \ref{discreteClifford2}, the computation of the inverse for $O_{x_j}'=\cosh\left(h\partial_{x_j}\right)$ is cumbersome and involves infinite sums or integral representations. However, in the case when periodic boundary conditions of the type $\underline{x}+h N{\bf v}_j=\underline{x}$ for certain $N \in \mathbb{N}$ are imposed on $h \BZ^n$ (see \cite{DHS96}, Section 5) it is possible to compute explicitly $(O_{x_j}')^{-1}$ as a finite sum involving powers of $T_{h \vv_j}=\exp\left(h\partial_{x_j}\right)$.

On the other hand, contrary to Example \ref{differenceDiracOperators}, the operators $x'$ and $D'$
obtained {\it viz} the following intertwining properties on $\mathcal{P}$:
\begin{eqnarray}
\label{intertwiningUmbral}\Psi_{\underline{x}} D=D'\Psi_{\underline{x}}, &  \Psi_{\underline{x}} x=x'\Psi_{\underline{x}}, & \Psi_{\underline{x}} E=E'\Psi_{\underline{x}}
\end{eqnarray}
concern with the nearest neighbor points together with all the points contained in each direction $h{\bf v}_j$.
Hereby, $\Psi_{\underline{x}}$ is the Sheffer map introduced in Section \ref{umbralC}.

Here we would also like to stress that $\Delta_{2h}=-(D')^2$ is supported on $(2h)\BZ^n$. So, the periodicity as well as the coarsening of lattice is the price that we must pay in order to get discrete Clifford analysis as a refinement of discrete harmonic analysis underlying the orthogonal group $O(n)$.

\subsection{Orthosymplectic Lie Algebra
Representation}\label{osp(1|2)}

The main objective of this subsection is to gather a fully description for the Clifford operators defined on Subsection \ref{umbralCliffordOperators} as a
representation of the orthosymplectic Lie algebra $\mathfrak{osp}(1|2)$. We will start to recall some basic definitions underlying the Lie algebra setting. A comprehensive survey of this topic can be found in \cite{Howe85,FSS00}.

The orthosymplectic Lie algebra of type $\mathfrak{osp}(1|2)$ is defined as 
$$\mbox{span}\left\{{\bf p}^-,{\bf p}^+,{\bf q}\right\}\oplus
\mbox{span}\left\{{\bf r}^-,{\bf r}^+\right\}$$ equipped with the standard graded
commutator $\left[ \cdot, \cdot\right]$ such that ${\bf p}^-,{\bf p}^+,{\bf q},{\bf r}^-$ and ${\bf r}^+$ satisfy the following standard commutation relations (see e.g.
\cite{FSS00}):
\begin{eqnarray*}
\left[{\bf q},{\bf p}^\pm \right]=\pm {\bf p}^\pm, & \left[{\bf p}^+,{\bf p}^-\right]=2{\bf q},\\
\left[{\bf q},{\bf r}^\pm \right]=\pm \frac{1}{2}{\bf p}^\pm, & \left[{\bf r}^+,{\bf r}^-\right]=\frac{1}{2}{\bf q}, \\
\left[{\bf p}^\pm,{\bf r}^\mp \right]=-{\bf r}^\pm, & \left[{\bf
r}^\pm,{\bf r}^\pm\right]=\pm \frac{1}{2}{\bf p}^\pm.
\end{eqnarray*}
Here we would like to point out that on above construction, the Lie algebra $\mathfrak{sl}_2(\BR)$ appears as a refinement of $\mathfrak{osp}(1|2)$ in the sense that the canonical generators ${\bf p}^-,{\bf p}^+,{\bf q}$ itself generate $\mathfrak{sl}_2(\BR)$. In particular $\mathfrak{sl}_2(\BR)$ corresponds to the even part of $\mathfrak{osp}(1|2)$.

In terms of the operators $x',D'$ and
$E'+\frac{n}{2}{\bf id}$ the subsequent lemma gives rise to an isomorphic characterization of $\mathfrak{osp}(1|2)$ . We leave the proof of the following lemma to Appendix \ref{LieSuperAlgebraAppendix}.

\begin{lemma}[See Appendix \ref{LieSuperAlgebraAppendix}]\label{LieSuperAlgebra}
The operators $x'$,$D'$ and $E'+\frac{n}{2}{\bf id}$ generate a
finite-dimensional \index{Lie Algebra}Lie algebra in
$\mbox{End}(\mathcal{P})$. The remaining commutation relations are
\small{\begin{eqnarray*}
\label{DH2}\left[ x', (x')^2\right]=0, & \left[ x', -\Delta'\right]=-2D', & \left[E'+\frac{n}{2}{\bf id},x'\right]=x'\\
\label{DH3} \left[ D', (x')^2\right]=-2x', & \left[ D', -\Delta'\right]=0, & \left[E'+\frac{n}{2}{\bf id},D'\right]=-D' \\
\label{DH4} \left[ (x')^2, -\Delta'\right]=4\left(E'+\frac{n}{2}{\bf
id}\right), & \left[ E'+\frac{n}{2}{\bf id},
-(x')^2\right]=-2(x')^2, & \left[ E'+\frac{n}{2}{\bf id},
-\Delta'\right]=2\Delta'
\end{eqnarray*}}
\end{lemma}

Furthermore, the standard commutation relations for $\mathfrak{osp}(1|2)$ are obtained by considering the following normalization:
\begin{eqnarray*}
\begin{array}{ccccc}
{\bf p}^-=-\frac{1}{2}\Delta', & {\bf p}^+=-\frac{1}{2}(x')^2, &
{\bf q}=\frac{1}{2}\left(E'+\frac{n}{2}{\bf id}\right), & {\bf r}^+=\frac{1}{2 \sqrt{2}}ix',~~{\bf r}^-=\frac{1}{2
\sqrt{2}}i D',
\end{array}
\end{eqnarray*}
and moreover, ${\bf p}^+=\frac{1}{2}(x')^2$, ${\bf p}^-=\frac{1}{2}\Delta$ and ${\bf
q}=\frac{1}{2}\left(E'+\frac{n}{2}{\bf id}\right)$ correspond to the canonical generators of $\mathfrak{sl}_2(\BR)$.

In brief, the above description establishes a parallel with the {\it continuum} versions of Clifford analysis (cf.~\cite{DSS}) and harmonic analysis (cf. \cite{HT92}) as representations of $\mathfrak{osp}(1|2)$ and $\mathfrak{sl}_2(\BR)$, respectively. This also establishes a link with the celebrated Wigner quantum systems introduced by Wigner in \cite{Wigner50} in the sense that the description of the Clifford-valued operator in terms of the symmetries of $\mathfrak{osp}(1|2)$ allows us to describe the motion of a particle confinining the quantum harmonic oscillator with $n$ degrees of freedom.

\section{Almansi type theorems in (discrete) Clifford analysis}\label{AlmansiTheorems}

\subsection{Main Result}

In this section we will derive an Almansi type theorem based on replacements of the operators $(x')^k$ by polynomial type operators $P_k(x')$ such that the mapping $(D')^kP_k(x'):\ker D' \rightarrow \ker D'$ is an isomorphism. For a sake a simplicity, we leave for Appendix \ref{AlmansiTheoremsAppendix} the proofs of the technical results regarding the proof of the main result.

We will start to pointing out the following definitions:

\begin{definition}
Let $\Omega$ be a domain in $\BR^n$ and $k\in \BN$.
 A function
$f: \Omega\longrightarrow \cl_{0, n}$ is \textsl{umbral
~polymonogenic} of degree $k$ if $(D')^kf(\underline{x})=0$ for all $\underline{x}\in \Omega$. For $k=1$ $f: \Omega\longrightarrow \cl_{0, n}$ is called
{\it umbral monogenic}.
\end{definition}

\begin{definition}
A domain $\Omega\subset\mathbf R^n$ is \textsl{starlike}
with center $0$ if for each $\underline{x}\in \Omega$ $t\underline{x}\in\Omega$ holds
for any $0\le t \le 1$.
\end{definition}

For each $k\in\BN$, if the right inverse of $(D')^kP_k(x')$ on the range of $\ker D'$ exists we will define it by $Q_k': \ker D' \rightarrow \ker D'$, i.e.
\begin{eqnarray*}
(D')^kP_k(x')(Q_k'f)=f, & \mbox{for all}~~f \in \ker D'.
\end{eqnarray*}

Thus, the Almansi theorem can be formulated as follows:

\begin{theorem}\label{th:01}
Let $\Omega$ be a starlike domain in $\BR^n$ with center $0$. If $f$
is a umbral polymonogenic function of degree $k$ in $\Omega$, then
there exist unique functions $f_0,f_1, \ldots,f_{k-2},f_{k-1}$, each one umbral
monogenic in $\Omega$ such that
 \begin{equation}\label{eq1.32}
 f(\underline{x})=P_0(x')f_0(\underline{x})+ P_1(x')f_1(\underline{x}) + \cdots + P_{k-1}(x')f_{k-1}(\underline{x}).
\end{equation}
 Moreover, the umbral monogenic functions
$f_0,f_1, \ldots,f_{k-2},f_{k-1}$ are given by the following formulas:

\begin{equation}\label{eq1.33}
\small{\begin{array}{rcr}
 f_{k-1}(\underline{x}) &=& Q_{k-1}' (D')^{k-1} f(\underline{x})
   \\
f_{k-2}(\underline{x})&=& Q_{k-2}' (D')^{k-2}({\bf id}-P_{k-1}(x')Q_{k-1}'
(D')^{k-1}) f(\underline{x})
    \\
& \vdots
 \\
f_1(\underline{x}) &=& Q_1D'({\bf id}-P_2(x')Q_2(D')^2)\cdots({\bf
id}-P_{k-1}(x')Q_{k-1} (D')^{k-1})f(\underline{x})
  \\
f_0(\underline{x}) &=& ({\bf id}-P_1(x') Q_1' D')({\bf
id}-P_2(x')Q_2'(D')^2)\cdots({\bf id}-P_{k-1}(x')Q_{k-1} (D')^{k-1})
f(\underline{x})\rlap.
\end{array}}
\end{equation}

Conversely the sum in  (\ref{eq1.32}) with $f_0,f_1, \ldots,f_{k-2},f_{k-1}$
umbral monogenic in $\Omega$, defines a umbral polymonogenic
function of degree $k$ in $\Omega$.
\end{theorem}

Before proving Theorem \ref{th:01}, we need some preliminary results.

\begin{lemma}[See Appendix \ref{AlmansiTheoremsAppendix}]\label{loweringLemma}
Let $\Omega$ be a starlike domain in $\BR^n$ with center $0$. For any Clifford-valued function $f(\underline{x})$ in $\Omega$, the following relations hold for each $s \in \BN$:
\begin{eqnarray}
\label{Dxsf}D'((x')^s f(\underline{x}))=-2(x')^{s-1}U_{s}'f(\underline{x})+(-1)^s (x')^s D'f(\underline{x}),
\end{eqnarray}
where $$U_s'=\left\{\begin{array}{ccc}
k~{\bf id}, & \mbox{if}~s=2k \\
E'+(\frac{n}{2}+k){\bf id}, & \mbox{if}~s=2k+1
\end{array}\right..
$$
\end{lemma}

From the above lemma, the next proposition naturally follows:

\begin{proposition}[See Appendix \ref{AlmansiTheoremsAppendix}]\label{mappingProperty}
The iterated umbral Dirac operator $(D')^k$ has the mapping property
$$ (D')^k: (x')^s \ker D' \rightarrow (x')^{s-k} \ker D'$$
for any $s \geq k$. Hereby, for each $f(\underline{x}) \in \ker D'$,
\begin{eqnarray}\label{D'kx's}
(D')^k\left((x')^sf(\underline{x})\right)=(-2)^{k} (x')^{s-k}U_{s-k+1}'\ldots U_{s-1}'U_{s}'f(\underline{x}),
\end{eqnarray}
where the operators $U'_j$ are defined in Lemma \ref{loweringLemma}.
\end{proposition}

Let $\Omega$ be a starlike domain with center $0$. For any $s>0$, we
define the operator  $I_s: C^1(\Omega, \cl_{0, n})\longrightarrow
C^1(\Omega, \cl_{0, n})$ by
\begin{equation}\label{eq1.11}
I_s f(\underline{x})=\int_0^1 f(t\underline{x}) t^{s-1} dt.
\end{equation}

In addition, we set $E_{s}=s{\bf id}+E.$ For $s=0$ we write $E$ instead of $E_{0}$.

\begin{lemma}[cf.~\cite{MR02}]\label{le:4.1}
 Let $\underline{x}\in \BR^n$ and $\Omega$ be a domain with
$\Omega\supset [0, \underline{x}]$.  If $s>0$ and $f\in C^1(\Omega,
\cl_{0,n})$, then
\begin{equation}\label{eq1.880} f(\underline{x})=I_s E_{s} f (\underline{x})= E_{s} I_s f(\underline{x}).\end{equation}
\end{lemma}

Sloppily speaking, the family of maps $I_s: C^1(\Omega,
\cl_{0,n}) \rightarrow C^1(\Omega,
\cl_{0,n})$ can be viewed as certain sort of right inverse for the operator $Dx=\sum_{j,k=1}^n\e_j\e_k \partial_{x_j}x_k$ in $\ker D$.
Indeed, if $f$ is monogenic from Lemma \ref{x&D} $$D(xf(\underline{x}))=x(Df(\underline{x}))+D(xf(\underline{x}))=-2Ef(\underline{x})-nf(\underline{x})=-2 E_{n/2}f(\underline{x})$$ holds whenever $O_{\underline{x}}=\partial_{\underline{x}}$.

Finally, from Lemma \ref{le:4.1}, $f(\underline{x})=-2 E_{n/2}\left(-\frac{1}{2}I_{n/2}f(\underline{x})\right)=-\frac{1}{2}D(xI_{n/2}f(\underline{x}))$, showing that $-\frac{1}{2}I_{n/2}$ is a right inverse for $Dx$ on the range $\ker D$.

On the other hand, when restricted to $\mathcal{P}=\oplus_{k=0}^\infty \mathcal{P}_k$, where each $f_k \in \mathcal{P}_k$ is a Clifford-valued homogeneous polynomial of degree $k$ ( i.e. $f_k(t \underline{x})=t^kf_k(\underline{x})$), the family of mappings $I_s$ satisfy the equation
\begin{eqnarray*}
I_s f_k(\underline{x})=\int_{0}^1 f_k(\underline{x})t^{k+s-1}dt=\frac{1}{k+s}f_k(\underline{x}), & \mbox{for all}~~f_k \in \mathcal{P}_k.
\end{eqnarray*}

Hence Lemma \ref{le:4.1}  remains true for $E'+s{\bf id}=\Psi_{\underline{x}}^{-1}E_s\Psi_{\underline{x}}$ in $\mathcal{P}$, since it holds componentwise.

\begin{lemma}[See Appendix \ref{AlmansiTheoremsAppendix}]\label{Is'Es'}
There exists $I_s': \mathcal{P} \rightarrow \mathcal{P}$ such that
$$(E'+s{\bf id})I_s'={\bf id}=I_s'(E'+s{\bf id}).$$
\end{lemma}

The next lemma will be also important on the sequel

\begin{lemma}[See Appendix \ref{AlmansiTheoremsAppendix}]\label{le:51}
  If $f\in \mathcal{P}$, then
\begin{equation}\label{eq3.1}
 D' I'_s f(\underline{x})=I'_{s+1} D' f(\underline{x}).
 \end{equation}
\end{lemma}

For any $k \in \BN_0$, denote by $Q_k'=\left(-\frac{1}{2}\right)^k(U_k')^{-1}(U_{k-1}')^{-1}\ldots (U_1')^{-1}$, where
\begin{eqnarray}
\label{Ts-1}(U_s')^{-1}=\left\{\begin{array}{ccc}
\frac{1}{k}~{\bf id}, & \mbox{if}~s=2k \\ \ \\
I'_{\frac{n}{2}+k}, & \mbox{if}~s=2k+1
\end{array}\right..
\end{eqnarray}

As direct consequence of (\ref{eq3.1}), we find that $I'_s f(\underline{x})$ is umbral monogenic whenever $f(\underline{x})$ is umbral monogenic. From the
definition of $Q_k'$ we thus obtain
 \begin{equation}\label{eq3.11}
 Q_k'(\ker D')=\ker D'.
 \end{equation}

Then the following lemma holds:

\begin{lemma}[See Appendix \ref{AlmansiTheoremsAppendix}]\label{MappingPropertyDkXk}
For any umbral monogenic
function $f$ in $\Omega$,
$$ (D')^k \left[(x')^k Q_k' f(\underline{x})\right]= f(\underline{x}), \quad \underline{x}\in \Omega.$$
\end{lemma}

Now we come to the proof of our main theorem for $P_k(x')=(x')^k$:

\proof[Proof of Theorem \ref{th:01}]
 It is
sufficient to show that
$$
\ker (D')^k=\ker (D')^{k-1}+P_{k-1}(x') \ker D', \quad k\in\BN,
$$
where $P_{k-1}(x')=(x')^{k-1}$. Notice that Lemma \ref{MappingPropertyDkXk} %
states that
\begin{equation}\label{eq3.3}
 (D')^k P_k(x') Q_k'={\bf id}.
\end{equation}

We divide the proof into two parts:
\begin{enumerate}
\item[(i)] $\ker (D')^k\supset \ker (D')^{k-1}+P_{k-1}(x') \ker D'$.
Since $\ker (D')^{k-1}\subset \ker (D')^{k}$, we need only to show $P_{k-1}(x') \ker D'\subset \ker (D')^k$. For any $g\in \ker D'$, by (\ref{eq3.3}) and
(\ref{eq3.11})  we have
\begin{eqnarray*}
 (D')^k (P_{k-1}(x')g) = D' ((D')^{k-1} P_{k-1}(x') Q_{k-1}' ) (Q_{k-1}')^{-1} g
= D' (Q_{k-1}')^{-1} g=0.
 \end{eqnarray*}
\item[(ii)]
$\ker (D')^k\subset \ker (D')^{k-1}+P_{k-1}(x') \ker D'$.

For any $f\in \ker (D')^k$, we have the decomposition
$$f=({\bf id}-P_{k-1}(x') Q_{k-1}' (D')^{k-1} )f + P_{k-1}(x') (Q_{k-1}' (D')^{k-1} f).$$
We will show that the first summand above is in $\ker (D')^{k-1}$ and the item in
the braces of the second summand is in $\ker D'$. This can be verified
directly. First,
\begin{eqnarray*}
 (D')^{k-1}({\bf id}-P_{k-1}(x') Q_{k-1} (D')^{k-1} )f=\\=
((D')^{k-1}- ((D')^{k-1} P_{k-1}(x') Q_{k-1}') (D')^{k-1} )f\\= ((D')^{k-1}-  (D')^{k-1} )f=0.
 \end{eqnarray*}
Next, since $(D')^{k-1}f\in \ker D'$ and $Q_{k-1}' \ker D'\subset \ker D'$, we
have $ Q_{k-1}' (D')^{k-1} f\in \ker D'$, as desired.
\end{enumerate}

\noindent This proves that $\ker (D')^k=\ker (D')^{k-1}+P_{k-1}(x') \ker D'$. By induction,
we can easily deduce that $\ker (D')^k=\ker D'+ P_1(x')\ker D' + P_2(x')\ker D' + \ldots + P_{k-1}(x')\ker D'$.

\medskip
Next we prove that for any $f\in \ker (D')^k$ the decomposition
$$ f= g+P_{k-1}(x')f_k, \quad g\in \ker (D')^{k-1}, f_k\in \ker D'$$
is unique. In fact, for such a decomposition, applying $(D')^{k-1}$
on both sides we obtain
\begin{eqnarray*}
(D')^{k-1} f &=& (D')^{k-1} g+(D')^{k-1} P_{k-1}(x') f_k
\\
&=& (D')^{k-1} P_{k-1}(x') Q_{k-1}' (Q_{k-1}')^{-1}  f_1
\\
&=& (Q_{k-1}')^{-1} f_k.
\end{eqnarray*}
Therefore
$f_k=Q_{k-1}'(D')^{k-1} f,$ so that
$$g=f-P_{k-1}(x')f_k=({\bf id}-P_{k-1}(x')Q_{k-1}'(D')^{k-1}) f.$$
Thus equations (\ref{eq1.32}) and (\ref{eq1.33}) follows by induction.

To prove the converse, we see from equation (\ref{D'kx's}) of Lemma \ref{mappingProperty} that $(D')^{k+1} (x')^k \ker D'=\{0\}$ holds for any
$k\in \BN$.

Replacing $k$ by $j$,  we have
$$
(D')^{k} (x')^j \ker D'=\{0\}$$ for any $k>j$.
 \qed

The proof of the above theorem can be interpreted as the following infinite triangle on which the subspaces $\ker (D')^k$ are despicted into columns. Each element of the triangle given by Proposition \ref{mappingProperty} corresponds to the action of $\mathfrak{osp}(1|2)\times O(n)$ on rows and columns:
$$\begin{array}{cccccccccccc}
\{ 0\}& &\ker D' &   &\ker (D')^2 & & \ker (D')^3 &  & \ker (D')^4  &   & \ldots
\\ \ \\
\{ 0\}& &\ker D' &\xrightarrow{x'}   &x'\ker D' & \xrightarrow{x'} & (x')^2\ker D' & \xrightarrow{x'} & (x')^3\ker D'    & & \ldots
\\
& &\downarrow{D'}& & \downarrow{D'} &  & \downarrow{D'} &  & \downarrow{D'} &  &   \\
& &\{0\}&   & \ker D' & \xrightarrow{x'} &  x'\ker D' & \xrightarrow{x'} & (x')^2\ker D'  & &  \ldots \\
& & &  &  \downarrow{D'} & &  \downarrow{D'} &  & \downarrow{D'} &  &  \\
& & &    & \{0 \} & & {\ker D'}  & \xrightarrow{x'} & x'\ker D'  &   & \ldots \\
& & &  &  &  & \downarrow{D'} &  & \downarrow{D'} &  &  \\
& & &    &  & & \{0\} &  & \ker D'  & &  \ldots \\
& & &  &  &  &  &  & \downarrow{D'} &  &  \\
& & &    &  & &  &  &  \{0\}  & & \ldots \\
& & &  &  &  &  &  &    \\
& & &    &  & &  &  &   & &\ldots \\
\end{array}$$
In those actions, the operator $x'$ shifts all the spaces in the same row to the right while the operator $D'$ shifts all the spaces in the same column down.
In particular, the $(k+1)-$line of the above diagram corresponds to the action of $(D')^k$ on the subspaces $(x')^s\ker D'$ represented in  $(s+1)-$column.

The next important step is the passage from the homogeneous operator of degree $k$, $(x')^k$, to a general polynomial type operator $P_k(x')$ with the mapping property $P_k(x'): \ker (D')^k \rightarrow \ker D'$. The corollary below gives a possible generalization for the construction of $P_k(x')$:

\begin{corollary}\label{th:01B}
If $P_k(x')=A_k'~(x')^k+ R_k(x')$ where $A_k'$ is a Hilbert-Schmidt operator acting on $\mathcal{P}$ that satisfy the graded commuting relation $[A_k',D']=a_kD'$ for some $a_k \in \BR$ and
\begin{eqnarray*}
(D')^k\left(R_k(x')f(\underline{x})\right)=0, & \mbox{for all}~~f \in \ker D',
\end{eqnarray*}
then Theorem \ref{th:01} fulfils whenever the eigenvalues of the operator $A_k'$ are greater than $k a_k$.
\end{corollary}

\proof
Starting from the definition of $P_k(x')$ and using induction on $k \in \BN_0$, the assumptions for $A_k'$ and $R_k(x')$ lead to
\begin{eqnarray*}
(D')^k(P_k(x')f(\underline{x}))=(D')^k(A_k'(x')^kf(\underline{x}))=(-ka_k{\bf id}+A_k')(D')^k\left((x')^kf(\underline{x})\right).
\end{eqnarray*}
whenever $f$ belongs to $\ker D'$ (i.e. $f$ is umbral monogenic).

From direct application of Proposition \ref{mappingProperty}, the later equation becomes then
$$(D')^k(P_k(x')f(\underline{x}))=(-2)^{k} (-ka_k{\bf id}+A_k')U_{1}'\ldots U_{k-1}'U_{k}'f(\underline{x})),$$
where the operators $U'_j$ are defined in Proposition \ref{loweringLemma}.

Replacement of $f(\underline{x})$ by $S_k'f(\underline{x})=\left(-\frac{1}{2}\right)^{k}(U_{k}')^{-1}(U_{k-1}')^{-1}\ldots (U_{0}')^{-1}f(\underline{x})$, on the above equation
results in $$(D')^k(P_k(x')S_k'f(\underline{x}))=(-ka_k{\bf id}+A_k')f(\underline{x}).$$
Hereby $(U_{s}')^{-1}$ are defined via equation (\ref{Ts-1}).

Now it remains to show that $-ka_k{\bf id}+A_k'$ is invertible ensuring that $Q_k'=S_k'(-ka_k{\bf id}+A_k')^{-1}$ is a right inverse for $(D')^kP_k(x'): \ker D' \rightarrow \ker D'$.

If $A_k'$ is a multiple of ${\bf id}$, $a_k=0$ and hence $A_k'$ is invertible and the proof of Corollary \ref{th:01B} follows from Theorem \ref{th:01}.

 Otherwise, since $A_k'$ is a Hilbert-Schmidt operator acting on $\mathcal{P}$ we conclude that $A_k'$ has discrete spectra. Then, analogously to the proof of Lemma \ref{Is'Es'} (see Appendix \ref{AlmansiTheoremsAppendix}) $A_k'$ is given by following the series expansion
$$ A_k' f(\underline{x})=\sum_{s=0}^\infty\lambda_{k,s}f_s(\underline{x}),$$
where $\lambda_{k,s}\in \BR$ correspond to the eigenvalues of $A_k'$.

Thus $-ka_k{\bf id}+A_k'$ is invertible whenever $-k a_k+\lambda_{k,s}$ is positive, that is $\lambda_{k,s}>k a_k$.

Finally, using the same order of ideas of the proof of Theorem \ref{th:01}, induction arguments lead to the following infinite triangle
$$\small{\begin{array}{cccccccccccc}
\{ 0\}& &\ker D' &   &\ker (D')^2 & & \ker (D')^3 &  & \ker (D')^4  &  &  \ldots
\\ \ \\
\{ 0\}& &P_0(x')\ker D' &\xrightarrow{x'}   &P_1(x')\ker D' & \xrightarrow{x'} & P_2(x')\ker D' & \xrightarrow{x'} & P_3(x')\ker D'  &  & \ldots
\\
& &\downarrow{D'}& & \downarrow{D'} &  & \downarrow{D'} &  & \downarrow{D'} &  &   \\
& &\{0\}&   & P_0(x')\ker D' & \xrightarrow{x'} &  P_1(x')\ker D' &  \xrightarrow{x'} & P_2(x')\ker D' & &\ldots \\
& & &  &  \downarrow{D'} & &  \downarrow{D'} &  & \downarrow{D'} &  & \\
& & &    & \{0 \} & & P_0(x'){\ker D'}  & \xrightarrow{x'} & P_1(x')\ker D'  &  & \ldots \\
& & &  &  &  & \downarrow{D'} &  & \downarrow{D'} &  &   \\
& & &    &  & & \{0\} &  & P_0(x')\ker D' & & \ldots     \\
& & &  &  &  &  &  & \downarrow{D'} &   \\
& & &    &  & &  &  &  \{0\}  &  & \ldots \\
& & &  &  &  &  &  &  &   \\
& & &    &  & &  &  &    & & \ldots \\
\end{array}}$$

This yields the following direct sum decomposition of $\ker (D')^k$:
\begin{eqnarray*}
\ker (D')^k&=&\ker (D')^{k-1}\oplus P_{k-1}(x')\ker D' \\
&=&\ker (D')^{k-2}\oplus P_{k-2}(x')\ker D'\oplus P_{k-1}(x')\ker D'\\
&=&\ldots \\
&=&P_0(x')\ker D' \oplus P_1(x')\ker D'\oplus \ldots \oplus P_{k-1}(x')\ker D'.
\end{eqnarray*}
concluding in this way the proof of Corollary \ref{th:01B}.
\qed

We will end this section by establishing a parallel between our approach and the approaches of \textsc{Ryan} (cf. \cite{Ryan90}), \textsc{Malonek \& Ren} (cf. \cite{MR02,MAGU}) and \textsc{Faustino \& K\"ahler} (cf. \cite{FK07}).

Recall that Fischer decomposition (\cite{DSS}, Theorem 1.10.1) states the spaces of homogeneous polynomials $\mathcal{P}_k$ are splitted in spherical monogenics pieces with degree not exceeding $k$:
$$ \mathcal{P}_k=\sum_{s=0}^{k} \bigoplus x^s \left(\mathcal{P}_{k-s} \cap \ker D \right).$$

Moreover, from the mapping property given by Lemma \ref{MappingPropertyDkXk} each $P_k \in \mathcal{P}_k$ belongs to $\ker D^{k+1}$ and hence from the intertwining property given by relations (\ref{intertwiningUmbral}) the Clifford-valued polynomial of degree $k$ given by $\Psi_{\underline{x}}P_k(\underline{x})=P_k(\underline{x}'){\bf 1}$ belongs to $\ker (D')^{k+1}$. Hence the following direct sum decomposition of $\ker (D')^{k+1}$: 
$$ \ker (D')^{k+1}=\ker D'\oplus x'\ker D' \oplus (x')^2\ker D' \oplus \ldots \oplus (x')^{k}\ker D',$$
comprise the approaches of \textsc{Ryan} (cf. \cite{Ryan90}), \textsc{Malonek \& Ren} (cf. \cite{MR02}) (i.e. for $\Psi_{\underline{x}}={\bf id}$) as well as the Fischer decomposition in terms of forward Dirac operators obtained by \textsc{Faustino \& K\"ahler} in \cite{FK07} if we consider the operators introduced in Example \ref{differenceDiracOperators}. The flexibility of this approach allows us also to get the Fischer decomposition for several classes of finite difference operators like the finite difference operators considered in Example \ref{discreteClifford2}.

\begin{remark}
The replacement of $(x')^k$ by the polynomial type operators $P_k(x')$ given by Corollary \ref{th:01B} gives a parallel in {\it continuum} with the decomposition in terms of iterated kernels obtained by \textsc{Ren \& Malonek} (cf.~\cite{MR02}) on which the operators $P_k(x')$ shall be interpreted as quantizations of Clifford-valued polynomials of degree $k$.
\end{remark}
\subsection{Parallelism with the Quantum Harmonic Oscillator}\label{HarmonicOscillator}

We will finish this section by turning out our attention to the quantum harmonic oscillator given by the following Hamiltonian written in terms of the potential operator $V_\hbar(x')=-\frac{1}{2}(x')^2-\frac{\hbar}{2}x'+\frac{\hbar^2}{8}\left( \Gamma'-\frac{n}{2}{\bf id}\right)$:
\begin{eqnarray*}
\mathcal{H}'_\hbar=-\frac{1}{2}\Delta'+V_\hbar(x'), & \mbox{with}~\hbar \in \BR.
\end{eqnarray*}
Hereby $\Gamma'=-x'D'-E'$ corresponds to the umbral counterpart of the spherical Dirac operator (cf. \cite{DSS}).

In order to analyze the $\mathfrak{sl}_2(\BR)$ symmetries of $\mathcal{H}'_\hbar$, we further introduce the following auxiliar operator acting on the algebra of Clifford-valued polynomials $\mathcal{P}$:
\begin{eqnarray*}
\mathcal{J}'_\hbar=\frac{\hbar}{4} D'+\frac{1}{2}\left(E'+\frac{n}{2}{\bf id}\right).
\end{eqnarray*}
 The subsequent proposition gives a
description of the Lie algebra symmetries underlying $\mathcal{H}'_\hbar$, showing that ${\bf p}^+=V_\hbar(x')$, ${\bf p}^-=-\frac{\Delta'}{2}$ and ${\bf q}=\frac{\hbar}{4}D'+E'+\frac{n}{2}{\bf id}$ correspond to the canonical generators of $\mathfrak{sl}_2(\BR).$

We start with the following lemma:
\begin{lemma}\label{commutingGamma}
When acting on $\mathcal{P}$, the operator $\Gamma'$ commute with the operators $E'$ and $\Delta'$:
\begin{eqnarray*}
[E',\Gamma']=0, & [\Delta',\Gamma']=0.
\end{eqnarray*}
\end{lemma}

\proof
For the proof of $[E',\Gamma']=0$ it remains to show that $[E',x'D']=0$ since from definition $[E',\Gamma']=[E',-x'D'-E']=-[E',x'D']$.

From Lemma \ref{LieSuperAlgebra} we get $[E',x']=x'$ and $[E',D']=-D'$. This leads to
$$ E'(x'D')=(x'+x'E')D'=x'D'+(-x'D'+x'D'E')=(x'D')E',$$
or equivalently $[E',x'D']=0$, as desired.

In order to show that $[\Gamma',\Delta']=0$, we recall the relations $[\Delta',x']=2D'$, $[\Delta',E']=2\Delta'$ and $[\Delta',D']=0$ that follow from Lemma \ref{LieSuperAlgebra}. This shows that
$$ \Delta'(x'D')=(2D'+x'\Delta')D'=-2\Delta'+(x'D')\Delta', $$
and hence, $$\Delta'(x'D'+E')=-2\Delta'+(x'D')\Delta'+2\Delta'+E'\Delta'=(x'D'+E')\Delta'.$$ Finally, taking into account the definition of $\Gamma'$ the above equation is equivalent to $[\Delta',\Gamma']=0,$ as desired.
\qed

\begin{lemma}\label{LieAlgebraVx'}
When acting on $\mathcal{P}$, the elements $\frac{\Delta'}{2}$,~$V_\hbar(x')$ and $\mathcal{J}_\hbar'$ are the canonical generators of the Lie algebra $\mathfrak{sl}_2(\BR)$.
The remaining commutation relations are
\begin{eqnarray*}
\left[\frac{\Delta'}{2},V_\hbar(x')\right]=\mathcal{J}'_\hbar, & \left[\mathcal{J}'_\hbar,V_\hbar(x')\right]=V_\hbar(x'),
& \left[\mathcal{J}'_\hbar,\frac{\Delta'}{2}\right]=\frac{\Delta'}{2}.
\end{eqnarray*}
\end{lemma}

\proof
Recall that from Lemma \ref{LieSuperAlgebra},
${\bf p}^-=-\frac{\Delta'}{2}$ ${\bf p}^+=\frac{(x')^2}{2}$ and ${\bf q}=\mathcal{J}'_\hbar$ are the canonical
 generators of $\mathfrak{sl}_2(\BR)$:
 \begin{eqnarray*}
\left[{\bf p}^-,{\bf p}^+ \right]={\bf q}, & \left[{\bf q},{\bf p}^-\right]=-{\bf p}^-, &
\left[{\bf q},{\bf p}^+
\right]={\bf q}.
 \end{eqnarray*}
Moreover $\left[x',{\bf p}^- \right]=-D'$, $\left[D',{\bf p}^+ \right]=x'$, $[{\bf q},x']=\frac{1}{2}x'$ and $\left[D',{\bf p}^- \right]=0=\left[x',{\bf p}^+ \right]$.
 Taking into account that
 \begin{center}
  $V_\hbar(x')=-{\bf p}^+-\frac{\hbar}{2}x'+\frac{\hbar^2}{8}\left( \Gamma'-\frac{n}{2}{\bf id}\right)$
and $\mathcal{J}'_\hbar={\bf q}+\frac{\hbar}{4}D'$,
\end{center}
combination of the above relations with Lemma \ref{commutingGamma} results in the following identities in terms graded commuting relations:
\begin{eqnarray*}
\begin{array}{lll}
\left[ -{\bf p}^-,V_\hbar(x')\right]&=&\left[ {\bf p}^-,{\bf p}^+\right]+\left[ {\bf p}^-,-\frac{\hbar}{2}x'+\frac{\hbar^2}{8}\left( \Gamma'-\frac{n}{2}{\bf id}\right)\right] \\
&=&{\bf q}+\frac{\hbar}{2}D'; \\ \nonumber \\
\left[\mathcal{J}'_\hbar ,V_\hbar(x')\right]&=&
-[{\bf q},{\bf p}^+]-\frac{\hbar}{2}\left[ {\bf q},x'\right]-\frac{\hbar}{4}\left[ D',{\bf p}^+\right]-(\frac{\hbar}{4})^2[D',x'] \\
&=&-{\bf p}^+-\frac{\hbar}{2}x'+\frac{\hbar^2}{8}\left( \Gamma'-\frac{n}{2}{\bf id} \right)\\
&=&V_\hbar(x');
 \\ \nonumber \\
\left[\mathcal{J}'_\hbar ,-{\bf p}^-\right]&=&\left[ {\bf q},{\bf p}^-\right]-\left[ \frac{\hbar}{4}D',{\bf p}^-\right]\\
&=&-{\bf p}^-.
\end{array}
\end{eqnarray*}

This proves Lemma \ref{LieAlgebraVx'}.
\qed

\begin{proposition}\label{H'J'}
The operators $\mathcal{J}'_\hbar,\mathcal{H}'_\hbar \in
\mbox{End}(\mathcal{P})$ are interrelated by the following intertwining property:
$$\mathcal{H}'_\hbar\exp\left(V_\hbar(x')\right)\exp\left(-\frac{\Delta'}{2}\right)=-\exp\left(V_\hbar(x')\right)\exp\left(-\frac{\Delta'}{2}\right)\mathcal{J}_\hbar'.$$
\end{proposition}

\proof
From Lemma \ref{LieAlgebraVx'}, the elements $\frac{\Delta'}{2}$,~$V_\hbar(x')$ and $\mathcal{J}_\hbar'$
correspond to the canonical generators of $\mathfrak{sl}_2(\BR)$.

From the above relations, it follows from induction over $k\in \BN$
that
\begin{eqnarray*}
\left[\frac{\Delta'}{2},V_\hbar(x')^k
\right]=k\mathcal{J}_\hbar'\left(V_\hbar(x')\right)^{k-1}, & \left[\mathcal{J}_\hbar',V_\hbar(x')^k
\right]=kV_\hbar(x')^{k-1}\left(V_\hbar(x')\right)^{k-1},
\end{eqnarray*}
leading to
\begin{eqnarray*}
\left[\frac{\Delta'}{2},\exp\left(V_\hbar(x')\right)
\right]=\mathcal{J}_\hbar'\exp\left(V_\hbar(x')\right), & \left[\mathcal{J}_\hbar',\exp\left(V_\hbar(x')\right)
\right]=V_\hbar(x')\exp\left(V_\hbar(x')\right).
\end{eqnarray*}

Combining the above relations we get
\begin{center}
$\left[\frac{\Delta'}{2}+\mathcal{J}_\hbar',\exp\left(V_\hbar(x')\right)\right]
=\left(\mathcal{J}_\hbar'+V_\hbar(x')\right)\exp\left(V_\hbar(x')\right).$
\end{center}
This is equivalent to $\left(\frac{\Delta'}{2}-V_\hbar(x') \right)\exp\left(V_\hbar(x')\right)=\exp\left(V_\hbar(x')\right)\left( \frac{\Delta'}{2}+\mathcal{J}_\hbar' \right)$.

Not it remains to show that $\left( \frac{\Delta'}{2}+\mathcal{J}_\hbar' \right)\exp\left(-\frac{\Delta'}{2}\right)=\exp\left(-\frac{\Delta'}{2}\right)\mathcal{J}_\hbar'.$

This statement is then immediate from the relation $$\left[\mathcal{J}_\hbar',\exp\left(-\frac{\Delta'}{2}\right)  \right]=\left(-\frac{\Delta'}{2}\right)\exp\left(-\frac{\Delta'}{2}\right).$$

 Therefore
\begin{eqnarray*}
 \mathcal{H}_\hbar\exp\left(V_\hbar(x')\right)\left(-\frac{\Delta'}{2}\right)
 &=&-\exp\left(V_\hbar(x')\right)\left(  \frac{\Delta'}{2}+\mathcal{J}'_\hbar\right)\exp\left(-\frac{\Delta'}{2}\right)\\
 &=&-\exp\left(V_\hbar(x')\right)\exp\left(-\frac{\Delta'}{2}\right)\mathcal{J}'_\hbar,
\end{eqnarray*}
as desired.
 \qed

We will finish this section by establishing a parallel with the recent approach of \textsc{Render} (cf. \cite{Render08}).

According to the definition of Fischer inner pair (\cite{Render08}, page 315) Corollary \ref{th:01B} shows that when $A'_k$ is a multiple of the identity operator, the pair $(P_{2k}(x'),(\Delta')^k)$ corresponds to a quantization of the Fischer inner pair that completely determines a (discrete) Almansi decomposition for polyharmonic functions (cf. \cite{Render08}, Proposition 20). Indeed, if for any umbral polyharmonic of function of degree $k$ on $\Omega$ (i.e. $(\Delta')^kf(\underline{x})=0$ holds on $\Omega$) we take $P_{2k}(x')=\left(2V_\hbar(x')\right)^k$ it is straightforward from Corollary \ref{th:01B} that the following decompositions holds:
$$f(\underline{x})=P_0(x')f_0(\underline{x})+ P_2(x')f_1(\underline{x}) + \cdots + P_{2k-2}(x')f_{k-1}(\underline{x})$$
where $f_0,f_1, \ldots,f_{k-2},f_{k-1}$ are umbral harmonic functions on $\Omega$ (i.e. $\Delta'f_j(\underline{x})=0$ holds for each $j=0,\ldots,k-1$ on $\Omega$ ).

Thus, it is also possible to obtain explicit formulae analogue to (\ref{eq1.32}) for umbral harmonic functions
$f_0,f_1, \ldots,f_{k-1}$ by considering the even powers of $D'$.

\begin{remark}\label{SphericalHarmonicsRemark}
Based on Lemma \ref{LieAlgebraVx'}, it is clear from the above construction that the functions $f_0,f_1, \ldots,f_{k-1}$ obtained from Corollary \ref{th:01B} are solutions of the coupled system of equations:
\begin{eqnarray*}
\label{sphericalHarmonics}\Delta' f_k=0, & \mathcal{J}_{\hbar}' f_k=(\frac{k}{2}+\frac{n}{4})f_k.
\end{eqnarray*}
It is clear that in the limit $\hbar \leftarrow 0$ the above coupled system of equations approximate the umbral counterpart of spherical harmonics. In addition, from Proposition \ref{H'J'} the composite action of $-\exp\left(V_\hbar(x')\right)\exp\left(-\frac{\Delta'}{2}\right)$ on each $f_k$ span the eigenfunctions of $\mathcal{H}'_\hbar$.

On the other hand straightforward computations combined Lemma \ref{LieSuperAlgebra}
leads to the following graded commuting property:
$$\left[\frac{1}{2}\left(E'+\frac{n}{2}{\bf id} \right),\exp\left( \frac{\hbar}{2}D'\right)\right]=-\frac{\hbar}{2}D'\exp\left( \frac{\hbar}{2}D'\right).$$
This yields the following intertwining property when restricted to the algebra $\mathcal{P}$: $$\frac{1}{2}\left(E'+\frac{n}{2}{\bf id} \right)\exp\left( -\frac{\hbar}{2}D'\right)=\exp\left(- \frac{\hbar}{2}D'\right)\mathcal{J}'_\hbar.$$
showing that $\exp\left( -\frac{\hbar}{2}D'\right)$ maps the umbral harmonic polynomials of degree $k$ onto umbral counterparts of spherical harmonics of degree $k$.

In this case the action $\exp\left( -\frac{\hbar}{2}D'\right)$ on $\mathcal{P}$ plays a similar role to the inversion of the Wick operator in Segal-Bargmann spaces underlying nilpotent Lie groups (cf. \cite{CnopsKisil99}).
\end{remark}

\section{Concluding Remarks and Open Problems}
In this paper we introduce an algebraic framework that can be seen as a comprised model for Clifford analysis underlying the orthogonal group $O(n)$. This makes it possible to construct the associated operators and polynomials in the discrete setting starting from the equations and their solutions in {\it continuum}.
The intertwining properties given by relations (\ref{intertwiningUmbral}) at the level of $\mbox{End}(\mathcal{P})$ gives us a meaningful interpretation of classical and discrete Clifford analysis as two quantal systems on which the Sheffer operator $\Psi_{\underline{x}}$ acts as a gauge transformation preserving the canonical relations between both systems.

This approach can be interpreted as a merge between radial algebra approach proposed by \textsc{Sommen} \cite{Sommen97} to define Clifford analysis with the quantum mechanical approach for umbral calculus described by \textsc{Dimakis, Hoissen \& Striker} (cf.~\cite{DHS96}) and \textsc{Levi, Tempesta \& Winternitz} (cf. \cite{LTW04}). Based on the recent approach of \textsc{Tempesta} (cf. \cite{Tempesta08}) we believe that this approach shall also be useful to construct polynomials in hypercomplex variables possessing the Appell set property. In this direction, the recent approaches of \textsc{Malonek \& Tomaz} (cf. \cite{MT08}) \textsc{De Ridder,De Schepper,~K\"ahler \& Sommen} (cf.\cite{RSKS10}) and \textsc{Bock, G\"urlebeck, L\'avi$\check{c}$ka \& V. Sou$\check{c}$ek} (cf. \cite{BGLS10}) are beyond to the Sheffer set property.

Here we would like to stress that contrary to the approaches of \textsc{Malonek \& Tomaz} and \textsc{Bock, G\"urlebeck, L\'avi$\check{c}$ka \& Sou$\check{c}$ek} on it is almost clear that the considered operators are generators of $\mathfrak{sl}_2(\BR)$ (or alternatively $\mathfrak{sl}_2(\BC)$) and $\mbox{osp}(1|2)$ while the Appell sets are invariant under the action of the orthogonal group $O(n)$, in the approach of \textsc{De Ridder,De Schepper,~K\"ahler \& Sommen} it was not yet realized for which group the Appell sets (or more generally, the Sheffer sets) are invariant.

Based on the recent paper of \textsc{Brackx, De Schepper,Eelbode \& Sou$\check{c}$ek} (cf.~\cite{BGLS10}) and the preprint of \textsc{Faustino \& K\"ahler} (cf.~\cite{FK08}), we conjecture the following:
\begin{center}
 `{\it All Hermitian operators represented in terms of $\mathfrak{sl}_2(\BR)$ and $\mathfrak{osp}(1|2)$ generators in {\it continuum} cannot be represented by $\mathfrak{sl}_2(\BR)$ and $\mathfrak{osp}(1|2)$ generators in discrete but instead by quantum deformations of it'}.
 \end{center}
For a nice motivation on this direction we refer to \cite{Howe89} (see Section 2) and also \cite{faustino09} (see Subsection 3.3) on which such gap was undertaken.

In the proof of Almansi decomposition (Theorem \ref{th:01} and Corollary \ref{th:01B}), the iterated (umbral) Dirac operators $(D')^k$ play a central role. In comparison with \cite{DSS,FK07,faustino09,RSKS10} we prove a similar result using the decomposition of the subspaces $\ker (D')^k$ based on resolutions of $\mathfrak{osp}(1|2) \times O(n)$ instead of considering {\it a-priori} a Fischer inner product.

With this framework, Theorem \ref{th:01} shows that the decomposition of $\ker (D')^k$ in terms of $\mathfrak{osp}(1|2) \times O(n)$ pieces yield the subspaces $(x')^s\ker D'$ for $s=0,1,\ldots,k-1$.
Moreover, the replacement of $(x')^k$ by a polynomial type operator $P_k(x')$ in Corollary \ref{th:01B} gives an alternative interpretation for decomposing kernel approach proposed by \textsc{Malonek \& Ren} (cf.~\cite{MAGU}) as well as refines the Fischer inner pair technique used by \textsc{Render} in \cite{Render08} to prove the Almansi decomposition in terms umbral polyharmonic functions.

As it was observed along this paper the resulting approach based on representation of the Lie algebra $\mathfrak{osp}(1|2)$ as a refinement $\mathfrak{sl}_2(\mathbb{R})$ has a core of applications in quantum mechanics that can further be consider to study special functions in Clifford analysis that belong to Segal-Bargmann spaces (see \cite{CFK10} and references therein). From the border view of physics, we have shown in Subsection \ref{HarmonicOscillator} that the approach obtained by \textsc{Render} shall be described using a diagonalization in terms of $\mathfrak{sl}_2(\BR)$. Indeed, Proposition \ref{H'J'} and Remark \ref{SphericalHarmonicsRemark} explains the parallel between the Almansi decomposition of the subspaces $\ker (\Delta')^k$ and the separation of variables of quantum harmonic oscillators (cf. \cite{Turbiner88,SmirnovTurbiner95,Zhang08}).

One may further bring this technique in the future to construct new families of Appell/Sheffer sets for hypercomplex variables as well as to study Schr\"odinger equations on grids. At this stage, new families of discrete Clifford-valued polynomials like e.g. hypercomplex generalizations of Kravchuk polynomials (cf. \cite{ML01}) should appear.

\section*{acknowledgments}

The major part of this work was developed in the \textit{Center of Research and Development in Mathematics and Applications} of University of Aveiro (Portugal). The authors would like also to thank to the members of the \textit{Complex and Hypercomplex analysis} research group for the kind hospitality and for the active-research environment during their stay.

\appendix

\section{Umbral Clifford Analysis}\label{LieSuperAlgebraAppendix}

\subsection{Proof of Lemma \ref{LieSuperAlgebra}}
\proof
 Notice that the relations $\left[ D', \Delta'~\right]=0=\left[
x', -(x')^2\right]$ are then fulfilled since $(x')^2$ and $-\Delta'$ commute with
all elements of $\mbox{End}(\mathcal{P})$ (first and second relations of Lemma \ref{x&D}).

The proof of $\left[E'+\frac{n}{2}{\bf id},x'\right]=x'$ {and}
$\left[E'+\frac{n}{2}{\bf id},D'\right]=-D'$ follow straightforward from the Weyl-Heisenberg character of operators $x_j'$ and $O_{x_j}$. Straightforward application
of the above relations naturally leads to
\begin{eqnarray*}
\left[ E'+\frac{n}{2}{\bf id}, (x')^2\right]&=&\left(x' + x'\left(E'+\frac{n}{2}{\bf id}\right)\right)x' - x'\left(-x' +\left(E'+\frac{n}{2}{\bf id}\right) x'\right)\\ &=&2(x')^2,
\end{eqnarray*}
\begin{eqnarray*}
\left[ E'+\frac{n}{2}{\bf id}, -\Delta' \right]&=&\left(-D' +
D'\left(E'+\frac{n}{2}{\bf id}\right)\right)D' - D'\left(D' +E'+\frac{n}{2}D'\right)\\
&=&2\Delta'.
\end{eqnarray*}

Furthermore the relations $\left[E'+\frac{n}{2}{\bf id},x'\right]=x'$,
$\left[E'+\frac{n}{2}{\bf id},D'\right]=-D'$ together with the third anti-commuting
relation of Lemma \ref{x&D} lead to
\begin{eqnarray*}
D'(x')^2-(x')^2D'&=&\left(-2\left(E'+\frac{n}{2}{\bf id}\right)-x'D'\right)x' -x'\left(-2\left(E'+\frac{n}{2}{\bf id}\right)-D'x'\right) \nonumber \\
&=&-2\left[E'+\frac{n}{2}{\bf id},x'\right] \nonumber \\
&=&-2x' \\ \ \\
-x'\Delta'+\Delta'x'&=&\left(-2\left(E'+\frac{n}{2}{\bf id}\right)-D'x'\right)D' -D'\left(-2\left(E'+\frac{n}{2}{\bf id}\right)-x'D'\right) \nonumber \\
&=&-2\left[E'+\frac{n}{2}{\bf id},D'\right] \nonumber \\
&=&2D'.
\end{eqnarray*}
Finally, the combination of the relations $[E'+\frac{n}{2}{\bf id},x']=x'$ and $[E'+\frac{n}{2}{\bf id},D']=-D'$ with the third anti-commuting relation of Lemma \ref{x&D} leads to
\begin{eqnarray*}
-\Delta'(x')^2&=&D'\left(-2x' +(x')^2D'\right)\\
&=&-2D'x'+\left(-2x'+(x')^2D'\right)D'\\
&=&-2\{x',D'\}-(x')^2\Delta'\\
&=&4\left(E'+\frac{n}{2}{\bf id}\right)-(x')^2\Delta'.
\end{eqnarray*}
\qed

\section{Almansi-type theorems in (discrete) Clifford analysis}\label{AlmansiTheoremsAppendix}

\subsection{Proof of Lemma \ref{loweringLemma}}
\proof
We use induction to prove (\ref{Dxsf}). Since
$\{x',D'\}=x'D'+D'x'=-2\left(E'+\frac{n}{2}{\bf id}\right)$ and $Dg(\underline{x})=0$, we have
\begin{equation}\label{eq3.29}
D'(x' g(\underline{x}))=-2 \left(E'+\frac{n}{2}{\bf id}\right) g(\underline{x}).
\end{equation}

Next we show  that, for any $\underline{x}\in\Omega$ and $k\in\BN$,

\begin{equation}\label{eq3.39}
\begin{array}{rcl} D' ((x')^{2k} g(\underline{x})) &=& -2k (x')^{2k-1} g(x);
\\
D' ((x')^{2k-1} g(\underline{x})) &=& -2 (x')^{2(k-1)}  \left(E'+(\frac{n}{2}+k-1){\bf id}\right) g(\underline{x}).
\end{array}
\end{equation}
 This can be checked by induction. Assuming that (\ref{eq3.39})
 holds for $k$. we shall now prove it also holds for $k+1$.
We now apply the operator  $x'D'+D'x'=-2\left(E'+\frac{n}{2}{\bf id}\right)$ to the function
$(x')^{2k}g(x)$:
\begin{equation}\label{eq3.191}
x'D'((x')^{2k}g(\underline{x}))+D'x'((x')^{2k}g(\underline{x}))=-2\left(E'+\frac{n}{2}{\bf id}\right)((x')^{2k}g(\underline{x})).
\end{equation}
By the hypothesis of induction, the first term in the left is equal
to $-2k(x')^{2k}g(\underline{x})$, while the left side equals
$-2(x')^{2k}\left(E'+(\frac{n}{2}+2k){\bf id}\right)g(\underline{x})$ due to the fact that
$$\left(E'+\left(\frac{n}{2}+s\right){\bf id}\right)x'-x'\left(E'+\left(\frac{n}{2}+s\right){\bf id}\right)=\left[E'+\frac{n}{2}{\bf id},x'\right]=x',$$
holds for all $s>0$.

As a result,
\begin{equation}\label{eq3.392}
\begin{array}{rcl} D'((x')^{2k+1}g(\underline{x}))&=&-x'D'(x'^{2k}g(\underline{x}))-2\left(E'+\frac{n}{2}{\bf id}\right)((x')^{2k}g(\underline{x}))
\\&=&2k(x)'^{2k}g(\underline{x})-2 (x')^{2k} \left(E'+\left(\frac{n}{2}+2k\right){\bf id}\right) g(\underline{x})
\\&=&-2 (x')^{2k} \left(E'+\left(\frac{n}{2}+k\right){\bf id}\right) g(\underline{x}).
\end{array}
\end{equation}
This proves the second equality of (\ref{eq3.39}). The first
equality of (\ref{eq3.39}) can be proved similarly. This proves the
identities (\ref{eq3.39}).
\qed

\subsection{Proof of Proposition \ref{mappingProperty}}
\proof
In order to prove the above mapping property, we will derive (\ref{D'kx's}) using induction over $k \in \BN$.
First notice that for $k=1$, relation (\ref{D'kx's}) automatically fulfils according to Lemma \ref{loweringLemma}.

Next we assume that (\ref{D'kx's}) holds for any $k-1$, with $k \in \BN$. Hence, the action of $D'$ on both sides of (\ref{D'kx's}) combined with Lemma \ref{loweringLemma} results in
\begin{eqnarray*}
(D')^k\left((x')^sf(\underline{x})\right)=\\=(-2)^{k+1} (x')^{s-k-1}U_{s-k}'U_{s-k+1}'\ldots U_{s-1}'U_{s}'f(\underline{x})+(-1)^{s-k}(x')^{s-k}D'g_{s-k}(\underline{x}),
\end{eqnarray*}
with $g_{s-k}(\underline{x})=(-2)^kU_{s-k}'U_{s-k+1}'\ldots U_{s-1}'U_{s}'f(\underline{x})$.

Now it remains to show that $D'g_{s-k}(\underline{x})=0$.
If $j$ is even, $U_j'=\frac{j}{2}{\bf id}$ and hence $[U_j',D']=0$. Otherwise, from the second relations of \ref{osp(1|2)} $\left[E'+\frac{n}{2},D'\right]=D'$ combined with the definition of $U'_j$ for $j$ odd results in $[U_j',D']=-D'$.

Thus, we have $[U_j',D']=-\frac{(-1)^j-1}{2}D'$ for each $j \in \BN$ and moreover for each $g \in \ker D'$ the action of $D'$ on $U_j'g(\underline{x})$ is equal to
$$D'(U_j'g(\underline{x}))=\frac{1-(-1)^j}{2}D'g(\underline{x})+U_j'(D'g(\underline{x}))=0,$$
that is, $D'U'_j(\ker D') \subset \ker D'$ for each $j \in \BN$.

Finally, recursive application of the above relation leads to
$$
D'U_{s-k}'U_{s-k+1}'\ldots U_{s-1}'U_{s}'~(\ker D') \subset \ker D'.$$
and this results in $D'g_{s-k}(\underline{x})=0,$ as desired.
\qed

\subsection{Proof of Lemma \ref{Is'Es'}}
\proof
Take $f_k(\underline{x}) \in \mathcal{P}$ such that $E'f_k(\underline{x})=kf_k(\underline{x})$ holds for each $k\in \BN_0$. Hence for $s>0$, the operator
$E'+s{\bf id}$ has only positive eigenvalues of the form $\lambda=k+s$ which shows that the inverse of $E'+s{\bf id}$ tactically exists.

 For any $f(\underline{x})=\sum_{k=0}^\infty f_k(\underline{x}) \in \mathcal{P}$ and $s>0$, define $I_s': \mathcal{P}\rightarrow \mathcal{P}$ as the operator given by the series expansion
$$ I_s'f(\underline{x})=\sum_{k=0}^\infty \frac{1}{k+s}f_k(\underline{x}).$$

Then we have $$(E'+s{\bf id})(I_s'f(\underline{x}))=\sum_{k=0}^\infty \frac{1}{k+s}\left((E'+s{\bf id})f_k(\underline{x})\right)=\sum_{k=0}^\infty f_k(\underline{x})=f(\underline{x})$$
Recalling the definition of $\Psi_{\underline{x}}$ in Section \ref{umbralC}, the Clifford-valued polynomial $\left(\Psi_{\underline{x}}^{-1}f_k\right)(\underline{x})$ is homogeneous of degree $k$ and hence
 $$I_s'f_k(\underline{x})=\Psi_{\underline{x}}\left(I_s\Psi_{\underline{x}}^{-1}f_k(\underline{x})\right)
=\Psi_{\underline{x}}\left(\frac{1}{k+s}\Psi_{\underline{x}}^{-1}f_k(\underline{x})\right)=\frac{1}{k+s}f_k(\underline{x}).$$
leads to
$$I_s'\left((E'+s{\bf id})f(\underline{x})\right)=\sum_{k=0}^\infty I_s'((k+s)f_k(\underline{x}))=\sum_{k=0}^\infty f_k(\underline{x})=f(\underline{x}).$$

This shows that $I_s'$ is an inverse for the operator $E'+s{\bf id}$, as desired.
\qed

\subsection{Proof of Lemma \ref{le:51}}
\proof
Starting from Lemma \ref{osp(1|2)}, we have $-D'=[E'+\frac{n}{2}{\bf id},D']$, or equivalently,
$$-D'=E'D'-D'E'=\left(E'+s{\bf id}\right)D'-D'\left(E'+s{\bf id}\right)$$
by adding and subtracting $(s-\frac{n}{2})D'$ on both sides of the first equation.
This is equivalent to $D'E_s'=E'_{s+1}D'$, where $s \mapsto E'_s=E'+s{\bf id}$.

Using the fact that
$E'_s=\left({I'}_s\right)^{-1}$, we end up with
$$D' I'_s =I'_{s+1}E_
{s+1}' D' I'_s=I'_{s+1}D'E'_sI'_s= I'_{s+1} D'.$$
\qed

\subsection{Proof of Lemma \ref{MappingPropertyDkXk}}
\proof
Denote $g(\underline{x})= Q_k' f(\underline{x})$. From  (\ref{eq3.11}) and the definition of $Q_k'$,
 $g$ is umbral monogenic in $\Omega$ and hence
from equation (\ref{D'kx's}) of Proposition \ref{mappingProperty}, we know that
\begin{equation*}
(D')^k((x')^kg(\underline{x}))=(-2)^{k} U_{1}'\ldots U_{k-1}'U_{k}'g(\underline{x}).
 \end{equation*}

Thus
$
(D')^k\left((x')^kQ_k'f(\underline{x})\right)=f(\underline{x})$ follows directly from the above induced formulas.
\qed

\end{document}